\documentclass[11pt]{article}
\usepackage{mathrsfs}
\usepackage{amsthm}
\usepackage{amssymb}
\usepackage{amsmath}
\usepackage{graphicx}
\usepackage{color}
\usepackage{amsfonts}
\usepackage{float}
\usepackage{cite}
\usepackage[text={140mm,210mm},left=35mm,vmarginratio=1:1]{geometry}
\newtheorem{theorem}{Theorem}[section]
\newtheorem{remark}{Remark}[section]
\newtheorem{lemma}[theorem]{Lemma}

\newtheorem{definition}[theorem]{Definition}

\newtheorem{corollary}[theorem]{Corollary}

\numberwithin{equation}{section}
\normalsize

\begin{document}
\title{\textbf{Hydrodynamics of the generalized $N$-urn Ehrenfest model}}

\author{Xiaofeng Xue \thanks{\textbf{E-mail}: xfxue@bjtu.edu.cn \textbf{Address}: School of Science, Beijing Jiaotong University, Beijing 100044, China.}\\ Beijing Jiaotong University}

\date{}
\maketitle

\noindent {\bf Abstract:} In this paper we are concerned with a generalized $N$-urn Ehrenfest model, where balls keeps independent random walks between $N$ boxes uniformly laid on $[0, 1]$. After a proper scaling of the transition rates function of the aforesaid random walk, we derive the hydrodynamic limit of the model, i.e., the law of large numbers which the empirical measure of the model follows, under an assumption where the initial number of balls in each box independently follows a Poisson distribution. We show that the empirical measure of the model converges weakly to a deterministic measure with density driven by an integral equation. Furthermore, we derive non-equilibrium fluctuation of the model, i.e, the central limit theorem from the above hydrodynamic limit. We show that the non-equilibrium fluctuation of the model is driven by a measure-valued time-inhomogeneous generalized O-U process. At last, we prove a large deviation principle from the hydrodynamic limit under an assumption where the transition rates function from $[0, 1]\times [0, 1]$ to $[0, +\infty)$ of the aforesaid random walk  is a product of two marginal functions from $[0, 1]$ to $[0, +\infty)$.

\quad

\noindent {\bf Keywords:} hydrodynamic limit, N-urn Ehrenfest model, non-equilibrium fluctuation, large deviation.

\section{Introduction}\label{section one}
In this paper, we are concerned with the generalized $N$-urn Ehrenfest model. The model is a continuous-time Markov process defined as follows. Assuming that there are $N$ boxes symbolled as $1,2,\ldots,N$. At $t=0$, there are some balls in every boxes according to some probability distribution. Let $\lambda\in C^{1,1}\left([0, 1]\times [0, 1]\right)$ be a positive function on $[0,1]\times [0,1]$ with continuous partial derivatives with respect to both coordinates. For each $1\leq i\leq N$, let $\lambda_N(i)=\sum_{j=1}^N\lambda(\frac{i}{N},\frac{j}{N})$.
Each ball in box $i$ is chosen independently at rate $\frac{1}{N}\lambda_N(i)$. When a ball $A$ in box $i$ is chosen, then for each box $j$, $A$ is put into $j$ with probability $\frac{\lambda(\frac{i}{N}, \frac{j}{N})}{\lambda_N(i)}$. Note that $A$ may be put into the box where it comes from.

The above process can be defined equivalently according to its generator. Let $X=\{0,1,2,\ldots\}^N$ and each $x\in X$ be written as
$x=\big(x(1), x(2), \ldots, x(N)\big)$, then there is a unique Markov process $\{X_t^N\}_{t\geq 0}$ on $X$ with generator $\mathcal{L}_N$ given by
\begin{equation}\label{equ 1.1 generator}
\mathcal{L}_Nf(x)=\sum_{i=1}^N\sum_{j=1}^N\frac{x(i)\lambda(\frac{i}{N}, \frac{j}{N})}{N}\big(f(x^{i,j})-f(x)\big)
\end{equation}
for any sufficiently smooth $f$ on $X$ and $x\in X$, where $x^{i,j}\in X$ satisfies $x^{i,j}=x$ when $i=j$ while
\[
x^{i,j}(l)=
\begin{cases}
x(l), & \text{~if~}l\neq i,j,\\
x(i)-1, & \text{~if~} l=i ,\\
x(j)+1, & \text{~if~} l=j
\end{cases}
\]
when $i\neq j$. According to the definition of $\mathcal{L}_N$, for the generalized $N$-urn Ehrenfest model, the number of balls in the $i$th box at moment $t$ is given by $X_t^N(i)$.

When $\lambda\equiv 1$, the model reduces to the standard $N$-Urn Ehrenfest model, where each ball is put into each box at rate $\frac{1}{N}$.

There is a long history for the research of the classic Ehrenfest model (see References \cite{Blom1989, Lathrop2016, Meerson2018, Palacios1994} and so on), where $N=2$, since this $2$-urn model describes the exchange of gaseous molecules between two containers. The case where $N=3$ is investigated in Reference \cite{Chen2017}, where an expected hitting time is given via electric networks. The main result in \cite{Chen2017} is extended in References \cite{Cheng2020} and  \cite{Song2018}, where analogue results of the expected hitting times are computed for the cases where $N\geq 4$. Furthermore, the distributions of hitting times for cases where $N\geq 4$ are given in \cite{Cheng2020}.

In this paper, we are concerned with hydrodynamics of the generalized $N$-urn Ehrenfest model, i.e., limit theorems of $\frac{1}{N}\sum_{i=1}^NX_t^N(i)\delta_{\frac{i}{N}}(dx)$ as $N\rightarrow +\infty$, where $\delta_a(dx)$ is Dirac measure concentrated on $a$. The theory of hydrodynamics describes the phenomenon that the microscopic density field of the stochastic model is dominated macroscopically by a PDE after properly space time scaling. There is a long history for the investigation of hydrodynamics of models such as exclusion processes (see Chapter 4 of \cite{kipnis+landim99}), zero-range processes (see \cite{De1984}), voter models (see \cite{PresuttiSpohn83}) and so on. For a detailed survey of the theory of hydrodynamics, see Reference \cite{kipnis+landim99}.

Hydrodynamics are generally investigated for models which are mass conserved, i.e., the total number of particles is invariant, such as exclusion processes and zero-range processes. Since the generalized $N$-urn Ehrenfest model is also mass conserved, it is natural to ask whether this model performs a hydrodynamic limit, which is our motivation to investigate the topic of this paper.

For an interacting particle system with a non-stationary initial distribution, central limit theorem from the hydrodynamic of the model is called `non-equilibrium fluctuation', which is a popular and important topic in recent years. Non-equilibrium fluctuations have already been widely discussed for models such as exclusion processes, voter models, reaction-diffusion models, binary contact path processes and so on. A common phenomenon of these models is that the non-equilibrium fluctuation is driven by some kind of measure-valued generalized Ornstein-Uhlenbeck process defined in \cite{Holley1978}. For mathematical details, see \cite{Erhard2020, PresuttiSpohn83, Jara2018, XueZhao2020+}. Inspired by this references, we study non-equilibrium fluctuations of $N$-urn Ehrenfest models as a second step of our investigation. We will show that the non-equilibrium fluctuation of our model is driven by a measure-valued time-inhomogeneous generalized O-U process. For a precise statement of our result, see Section \ref{section two}.

The study of large deviations from hydrodynamics of interacting particle systems dates back to 1980s, when large deviations of simple exclusion processes are discussed in \cite{Kipnis1989}. Proofs of main results in \cite{Kipnis1989} given by Kipnis, Olla and Varadhan provide a later routine strategy for the research of large deviations from hydrodynamics and motivate many follow-up works such as \cite{Funaki2001, Gao2003, Lasinion1993} and so on. In this paper, inspired by \cite{Kipnis1989}, we prove a large deviation principle of our $N$-urn Ehrenfest model under an assumption that the transition rates function $\lambda(x, y)$ is a product of two marginal one-dimensional functions. For mathematical details, see Section \ref{section two}.

\section{Main results}\label{section two}
In this section we give our main results. For later use, we first introduce some notations and definitions. Throughout this paper, we let $T>0$ be a given moment. We use $\mathcal{S}$ to denote the dual of $C[0, 1]$ with the topology such that $\mathscr{A}_n\rightarrow \mathscr{A}$ in $\mathcal{S}$ if and only if
\[
\lim_{n\rightarrow +\infty}\mathscr{A}_n(f)=\mathscr{A}(f)
\]
for any $f\in C[0, 1]$. For a bounded signed measure $\mathfrak{m}$, we use $\langle\mathfrak{m}, f\rangle$ to denote $\int_{[0, 1]}f(x)\mathfrak{m}(dx)$ for any $f\in C[0, 1]$. Consequently, $\mathfrak{m}$ can be identified with an element in $\mathcal{S}$ such that $\mathfrak{m}(f)=\langle\mathfrak{m}, f\rangle$ for $f\in C[0, 1]$.

We use $D\big([0, T], \mathcal{S}\big)$ to denote the Skorokhod space of c\`{a}dl\`{a}g functions $F:[0, T]\rightarrow \mathcal{S}$. That is to say, for any $F\in D\big([0, T], \mathcal{S}\big)$ and $0<t<T$, $\lim_{s\downarrow t}F_s=F_t$ while $\lim_{s\uparrow t}F_s=F_{t-}$ exists in the sense that
\[
\lim_{s\downarrow t}F_s(f)=F_t(f) \text{~and~} \lim_{s\uparrow t}F_s(f)=F_{t-}(f)
\]
for any $f\in C[0,1]$.

For $0\leq t\leq T$ and $N\geq 1$, we define $\mu_t^N=\frac{1}{N}\sum_{i=1}^NX_t^N(i)\delta_{\frac{i}{N}}(dx)$, i.e.,
\[
\mu_t^N(f)=\langle\mu_t^N, f\rangle=\frac{1}{N}\sum_{i=1}^NX_t^N(i)f(\frac{i}{n})
\]
for any $f\in C[0,1]$. Then, $\mu^N:=\big\{\mu_t^N\big\}_{0\leq t\leq T}$ is a random element in $D\big([0, T], \mathcal{S}\big)$. We define $P_1, P_2: C[0,1]\rightarrow C[0,1]$ as linear operators such that
\[
(P_1f)(x)=\int_0^1\lambda(x,y)f(y)dy \text{\quad and \quad}(P_2f)(x)=f(x)\int_0^1\lambda(x,y)dy
\]
for any $f\in C[0,1]$ and $x\in [0,1]$. We define $P_1^{*}, P_2^{*}: \mathcal{S}\rightarrow \mathcal{S}$ as linear operators such that
\[
(P_1^{*}\mathscr{A})(f)=\mathscr{A}(P_1f) \text{\quad and \quad} (P_2^{*}\mathscr{A})(f)=\mathscr{A}(P_2f)
\]
for any $\mathscr{A}\in\mathcal{S}$ and $f\in C[0,1]$.

To give our main result of $\mu^N$, we need the following assumption.

\quad

\textbf{Assumption} (A): For all $N\geq 1$, $\{X_0^N(i)\}_{i=1}^N$ are independent and $X_0^N(i)$ follows a Poisson distribution with mean $\phi\left(\frac{i}{N}\right)$ for each $i$, where $\phi\in C[0,1]$.

It is easy to check that, under Assumption (A), $\frac{1}{N}\sum_{i=1}^NX_0^N(i)f\left(\frac{i}{N}\right)$ converges weakly to $\int_0^1 f(x)\phi(x)dx$ for any $f\in C[0,1]$ as $N$ grows to infinity.

\begin{definition}\label{definition 2.1}
$\mu=\{\mu_t\}_{0\leq t\leq T}\in D\left([0, T], \mathcal{S}\right)$ is called as a weak solution of the following Equation
\begin{equation}\label{equ 2.1 ODE on banach}
\begin{cases}
&\frac{d}{dt}\mu_t=\big(P_1^{*}-P_2^{*}\big)\mu_t, \text{\quad}0\leq t\leq T,\\
&\mu_0(dx)=\phi(x)dx,
\end{cases}
\end{equation}
if for any $0\leq t\leq T$ and $f\in C[0, 1]$,
\begin{equation}\label{equation 2.1 Mu with f}
\mu_t\left(f\right)=\int_0^1f(x)\phi(x)dx+\int_0^t \mu_s\left((P_1-P_2)f\right)ds.
\end{equation}
\end{definition}
Equation \eqref{equ 2.1 ODE on banach} has the following property.

\begin{lemma}\label{Existen and uniqueness}
There exists an unique weak solution $\mu\in D\left([0, T], \mathcal{S}\right)$ of Equation \eqref{equ 2.1 ODE on banach}. Furthermore, $\mu_t(dx)=\rho(t,x)dx$ for any $0\leq t\leq T$, where
\begin{equation}\label{equ 2.1 Integration Equation for density}
\begin{cases}
& \frac{d}{dt}\rho(t,x)=\int_0^1\lambda(y,x)\rho(t,y)dy-\rho(t,x)\int_0^1\lambda(x,y)dy, \text{\quad}0\leq t\leq T,\\
&\rho(0,x)=\phi(x).
\end{cases}
\end{equation}
\end{lemma}

We put the proof of Lemma \ref{Existen and uniqueness} in the appendix, which utilizes classic theory of ordinary differential equations on Banach spaces. Now we can state our first main result, which is the following limit theorem of $\{\mu^N\}_{N\geq 1}$.

\begin{theorem}\label{theorem 2.1 main result hydrolimit}
Under Assumption (A), $\mu^N$ converges weakly to the weak solution $\mu$ of Equation \eqref{equ 2.1 ODE on banach} as $N$ grows to infinity.
\end{theorem}

\quad

Our next result is regrading the non-equilibrium fluctuation, i.e., central limit theorem, corresponding to the hydrodynamic limit given in Theorem \ref{theorem 2.1 main result hydrolimit}.

To give our result, we first introduce some notations. For any $t\geq 0$ and $N\geq 1$, we use $V_t^N$ to denote $\frac{1}{\sqrt{N}}\sum_{i=1}^N\big(X_t^N(i)-{\rm E}\left(X_t^N(i)\right)\big)\delta_{\frac{i}{N}}(dx)$. Then, $V^N:=\{V_t^N\}_{0\leq t\leq T}$ is also a random element in  $D\big([0, T], \mathcal{S}\big)$. We use $\mathcal{W}=\{\mathcal{W}_t\}_{t\geq 0}$ to denote the $\left(C\left([0,1]\times[0,1]\right)\right)^\prime$-valued standard Brownian motion, i.e., for any $f\in C\left([0,1]\times[0,1]\right)$, $\{\mathcal{W}_t(f)\}_{t\geq 0}$ is a real-valued Brownian motion such that $\mathcal{W}_0(f)=0$ and
\[
{\rm Cov}\big(\mathcal{W}_t(f), \mathcal{W}_t(f)\big)=t\int_0^1\int_0^1f^2(x, y)dxdy
\]
for any $t\geq 0$ (see page 608 of \cite{Gon2014}).

For any $t\in [0, T]$, we define $b_t: C[0,1]\rightarrow C([0,1]\times[0,1])$ as the linear operator such that
\[
(b_tf)(x, y)=\sqrt{\rho(t,x)\lambda(x,y)}\big(f(y)-f(x)\big)
\]
for any $f\in C[0,1]$ and $x, y\in [0,1]$, where $\rho$ is defined as in Equation \eqref{equ 2.1 Integration Equation for density}. Then, for any $t\in [0, T]$, we define $b_t^*: \left(C\left([0,1]\times[0,1]\right)\right)^\prime\rightarrow \mathcal{S}$ as the linear operator such that
\[
\left(b_t^*\mu\right)(f)=\mu(b_tf)
\]
for any $f\in C[0,1]$ and $\mu\in \left(C\left([0,1]\times[0,1]\right)\right)^\prime$.

We also need to introduce the definition of a time-inhomogeneous generalized Ornstein-Uhlenbeck process. According to an analysis similar with that given in the proof of Theorem (1.4) of \cite{Holley1978}, there exists an unique stochastic element $V=\{V_t\}_{0\leq t\leq T}$ in $D\left([0, T], \mathcal{S}\right)$ satisfying:

(1) $\{V_t(H)\}_{0\leq t\leq T}$ is a real-valued continuous function for any $H\in C[0, 1]$.

(2) For any $H\in C[0,1]$ and any $G\in C_c^{\infty}(\mathbb{R})$,
\begin{align*}
\Bigg\{G\left(V_t(H)\right)-G\left(V_0(H)\right)&-\int_0^tG^\prime\left(V_s(H)\right)V_s\left((P_1-P_2)H\right)ds \\
&-\frac{1}{2}\int_0^tG^{\prime\prime}\left(V_s(H)\right)\|b_sH\|_2^2ds\Bigg\}_{0\leq t\leq T}
\end{align*}
is a martingale, where $\|f\|_2=\sqrt{\int_0^1\int_0^1 f^2(x, y)dxdy}$ for any $f\in L^2\left([0, 1]\times[0, 1]\right)$.

(3) $V_0(H) \text{~follows~}\mathbb{N}\left(0, \int_0^1H^2(x)\phi(x) dx\right)\text{~for any~}H\in C[0, 1]$, where $\mathbb{N}(\mu ,\sigma^2)$ is the normal distribution with mean $\mu$ and variance $\sigma^2$.

Then, it is natural to define the above $V$ as the solution to the $\mathcal{S}$-valued Equation:
\begin{equation}\label{equ 2.1}
\begin{cases}
&dV_t=(P_1^{*}-P_2^{*})V_tdt+b_t^*d\mathcal{W}_t, \text{\quad}0\leq t\leq T,\\
&V_0(H) \text{~follows~}\mathbb{N}\left(0, \int_0^1H^2(x)\phi(x) dx\right)\text{~for any~}H\in C[0, 1],\\
&V_0 \text{~is independent of~}\{\mathcal{W}_t\}_{t\geq 0},
\end{cases}
\end{equation}
where $\{\mathcal{W}_t\}_{t\geq 0}$ is the $\left(C\left([0,1]\times[0,1]\right)\right)^\prime$-valued Brownian motion as we have mentioned above.

Now we give our second main result, which gives non-equilibrium fluctuations of generalized $N$-urn Ehrenfest models.

\begin{theorem}\label{theorem 2.2 main result CLT}
Under Assumption (A), $V^N$ converges weakly to $V$ as $N\rightarrow +\infty$, where $V=\{V_t\}_{0\leq t\leq T}$ is the unique solution to Equation \eqref{equ 2.1}.
\end{theorem}

\begin{remark}\label{remark Solving OU}
The solution of Equation \eqref{equ 2.1} can be defined equivalently through solving this Equation directly. As we have mentioned in the proof of Lemma \ref{Existen and uniqueness}, it is reasonable to define
\[
e^{u(P_1-P_2)}=\sum_{k=0}^{+\infty}\frac{u^k\left(P_1-P_2\right)^k}{k!}
\]
and hence define $e^{u\left(P_1^{*}-P_2^{*}\right)}=\sum_{k=0}^{+\infty}\frac{u^k\left(P_1^{*}-P_2^{*}\right)^k}{k!}$
for any $u\in \mathbb{R}$. Let $V$ be the solution of Equation \eqref{equ 2.1}, then
\[
d\left(e^{-t\left(P_1^{*}-P_2^{*}\right)}V_t\right)=e^{-t\left(P_1^{*}-P_2^{*}\right)}b_t^{*}d\mathcal{W}_t
\]
and hence
\[
V_t=e^{t\left(P_1^{*}-P_2^{*}\right)}V_0+\int_0^t e^{(t-s)\left(P_1^{*}-P_2^{*}\right)}b_s^{*}d\mathcal{W}_s,
\]
i.e.,
\[
V_t(H)=V_0\left(e^{t(P_1-P_2)}H\right)+\eta_t(H)
\]
for any $H\in C[0, 1]$, where $\{\eta_t(H)\}_{0\leq t\leq T}$ is a continuous martingale independent of $V_0$ and satisfies
\[
\langle \eta(H), \eta(H)\rangle_t=\int_0^t \left\|b_s\left(e^{(t-s)(P_1-P_2)}(H)\right)\right\|_2^2ds,
\]
where $\|b_s(f)\|_2^2=\int_0^1\int_0^1 \rho(s,x)\lambda(x, y)\left(f(x)-f(y)\right)^2dxdy$ as we have defined above.
\end{remark}

According to Remark \ref{remark Solving OU}, we have the following result as a direct corollary of Theorem \ref{theorem 2.2 main result CLT}.
\begin{corollary}\label{Corollary 2.3}
Under Assumption (A), for any $t>0$ and $H\in C[0,1]$,
\[
V_t^N(H)=\frac{1}{\sqrt{N}}\sum_{i=1}^N\left(X_t^N(i)-{\rm E}\left(X_t^N(i)\right)\right)H(\frac{i}{N})
\]
converges weakly to $\mathbb{N}\big(0, \theta^2(t, H)\big)$ as $N\rightarrow+\infty$, where
\[
\theta(t,H)=\sqrt{\int_0^1 \left(e^{t(P_1-P_2)}H\right)^2(x)\phi(x)dx+\int_0^t\left\|b_s\left(e^{(t-s)(P_1-P_2)}(H)\right)\right\|_2^2ds}.
\]
\end{corollary}

Our third main result is regarding the large deviation principle from the hydrodynamic limit given in Theorem \ref{theorem 2.1 main result hydrolimit}. We first introduce some notations and definitions. For any $\nu\in \mathcal{S}$, we define
\[
I_{ini}(\nu)=\sup\left\{\nu(H)-\int_0^1\phi(x)\left(e^{H(x)}-1\right)dx:~H\in C[0, 1]\right\}.
\]
We define nonlinear operator $\mathcal{B}: C[0, 1]\rightarrow C[0 ,1]$ as
\[
\mathcal{B}f(x)=\int_0^1\lambda(x, y)\left(e^{f(y)-f(x)}-1\right)dy
\]
for any $f\in C[0, 1]$ and $x\in [0, 1]$.

For any $\mu=\{\mu_t\}_{0\leq t\leq T}\in D\big([0, T], \mathcal{S}\big)$, we define
\begin{align*}
I_{dyn}(\mu)=\sup\Bigg\{&\mu_T(G_T)-\mu_0(G_0)\\
&-\int_0^T\mu_s\left((\partial_s+\mathcal{B})G_s\right)ds:~G\in C^{1,0}\left([0, T]\times [0,1]\right)\Bigg\}.
\end{align*}

We use $D_0$ to denote the subset of $D\big([0, T], \mathcal{S}\big)$ consists of $\mu$ satisfying:

(1) There exists nonnegative $\psi\in C^{1, 0}\left([0, T]\times [0, 1]\right)$ such that $\psi(t,\cdot)(x)=\frac{d\mu_t}{dx}$ for $0\leq t\leq T$.

(2) $\int_0^1\partial_s\psi(s, x)dx=0$ for $0\leq s\leq T$.

For some technical reason, currently we need the following assumption to give our large deviation principle.

\textbf{Assumption} (B): There exists $\lambda_1, \lambda_2\in C[0, 1]$ such that $\lambda(x, y)=\lambda_1(x)\lambda_2(y)$ for any $x, y\in C[0, 1]$.

Now we give our main result.

\begin{theorem}\label{theorem 2.3 LDP}
Under Assumptions (A) and (B), for any closed set $C\subseteq D\big([0, T], \mathcal{S}\big)$,
\begin{equation}\label{equ LDP upper}
\limsup_{N\rightarrow+\infty}\frac{1}{N}\log P\left(\mu^N\in C\right)\leq -\inf_{\mu\in C}\left(I_{ini}(\mu_0)+I_{dyn}(\mu)\right)
\end{equation}
while for any open set $O\subseteq D\big([0, T], \mathcal{S}\big)$,
\begin{equation}\label{equ LDP lower}
\liminf_{N\rightarrow+\infty}\frac{1}{N}\log P\left(\mu^N\in O\right)\geq -\inf_{\mu\in D_0\cap O}\left(I_{ini}(\mu_0)+I_{dyn}(\mu)\right).
\end{equation}
\end{theorem}

\begin{remark}
(1)According to our current approach utilized in the proof of Theorem \ref{theorem 2.3 LDP}, Assumption (B) is required to prove the fact that $\mu\in D_0$ implies that there exists $G\in C\left([0, T]\times[0, 1]\right)$ such that
\begin{equation}\label{equ 2.7}
\partial_s\psi(s, x)=\int_0^1\psi(s, y)\lambda(y, x)e^{G_s(x)-G_s(y)}-\psi(s, x)\lambda(x, y)e^{G_s(y)-G_s(x)}dy
\end{equation}
for any $(s, x)\in[0, T]\times [0, 1]$, where $\mu_t(dx)=\psi(t, x)dx$ for $0\leq t\leq T$. We guess that $G$ defined in Equation \eqref{equ 2.7} still exists without Assumption (B) and consequently Theorem \ref{theorem 2.3 LDP} holds for general $\lambda\in C\left([0, 1]\times [0, 1]\right)$, We will work on this problem as a further investigation.

(2) We guess that $I_{dyn}(\mu)<+\infty$ if and only if $\mu\in D_0$. If this property holds, then Equation \eqref{equ LDP lower} can be improved to
\[
\liminf_{N\rightarrow+\infty}\frac{1}{N}\log P\left(\mu^N\in O\right)\geq -\inf_{\mu\in O}\left(I_{ini}(\mu_0)+I_{dyn}(\mu)\right).
\]
We will work on this problem as a further investigation.
\end{remark}

The proofs of Theorems \ref{theorem 2.1 main result hydrolimit} and \ref{theorem 2.2 main result CLT} are given in Sections \ref{section three} and \ref{section four} respectively. The strategies of the two proofs are similar. First we will show that $\{\mu^N\}_{N\geq 1}$ and $\{V^N\}_{N\geq 1}$ are tight. Then, by Dynkin’s martingale formula, we show that any limit point $\widetilde{\mu}$ of subsequence of $\{\mu^N\}_{N\geq 1}$ satisfies Equation \eqref{equation 2.1 Mu with f} for any $f\in C[0,1]$ while any limit point $\widetilde{V}$ of subsequence of $\{V^N\}_{N\geq 1}$ is a solution of Equation \eqref{equ 2.1}. Consequently, $\widetilde{V}=V$ (resp. $\widetilde{\mu}=\mu$) holds according to the uniqueness of the solution to Equation \eqref{equ 2.1} (resp. Equation \eqref{equ 2.1 ODE on banach}).

The proof of Theorem \ref{theorem 2.3 LDP} is divided into Sections \ref{section five} and \ref{section six}. We follows a strategy similar with that introduced in \cite{Kipnis1989}, where an exponential martingale plays the crucial role. The proof of Equation \eqref{equ LDP upper} is given in Section \ref{section five}, where we first utilize the aforementioned exponential martingale, Chebyshev's inequality and a minimax theorem given in \cite{Sion1958} to show that Equation \eqref{equ LDP upper} holds for all compact $K\subseteq D\big([0, T], \mathcal{S}\big)$. Then, we utilize the criteria introduced in \cite{Puhalskii1994} to show that $\{\mu^N\}_{N\geq 1}$ are exponential tight and consequently complete the proof. The proof of Equation \eqref{equ LDP lower} is given in Section \ref{section six}, where a critical step is the utilization of a generalized version of Girsanov's theorem given in \cite{Schuppen1974} to derive the law of large numbers $\{\mu^N\}_{N\geq 1}$ obeys under a transformed measure with the aforementioned exponential martingale as the R-N derivative with respect to the original measure of $\{X_t^N\}_{0\leq t\leq N}$. The above strategy has also been utilized in the study of moderate deviations of stochastic models, see References \cite{Gao2003, Xue2019+, XueZhao2020+b} as examples.

\section{Proof of Theorem \ref{theorem 2.1 main result hydrolimit}}\label{section three}
In this section we prove Theorem \ref{theorem 2.1 main result hydrolimit}. For simplicity, in proofs we use $o(1)$ to denote a deterministic sequence which converges to $0$ as $N\rightarrow+\infty$ while use $o_p(1)$ to denote a stochastic sequence which converges weakly to $0$ as $N\rightarrow+\infty$. We use $O(1)$ to denote a deterministic bounded sequence. First we show the tightness of $\{\mu^N\}_{N\geq 1}$.

\begin{lemma}\label{lemma 3.1}
Under Assumption (A), $\{\mu^N\}_{N\geq 1}$ is tight.
\end{lemma}

\proof

According to Theorem 4.1 of Reference \cite{Mitoma1983} and Aldous' criteria, to show the tightness of $\{\mu^N\}_{N\geq 1}$, we only need to check that the following two conditions hold.

(1) $\lim_{M\rightarrow+\infty}\limsup_{N\rightarrow+\infty}P\big(|\mu_t^N(f)|\geq M\big)=0$ for any $t\geq 0$ and $f\in C[0,1]$.

(2) For any $\epsilon>0$ and $f\in C[0, 1]$,
\[
\lim_{\delta\rightarrow 0}\limsup_{N\rightarrow+\infty}\sup_{\tau\in \mathcal{T}, s\leq\delta}P\big(|\mu_{\tau+s}^N(f)-\mu_\tau^N(f)|>\epsilon\big)=0,
\]
where $\mathcal{T}$ is the set of stopping times of $\{X_t^N\}_{t\geq 0}$ bounded by $T$.

To check Condition (1), we have
\begin{align*}
P\big(|\mu_t^N(f)|\geq M\big)\leq &P\big(|\mu_t^N(f)|\geq M, \frac{1}{N}\sum_{i=1}^NX_0^N(i)\leq 2\int_0^1 \phi(x)dx\big)\\
&+P\big(\frac{1}{N}\sum_{i=1}^NX_0^N(i)>2\int_0^1 \phi(x)dx\big).
\end{align*}
Let $\|f\|_{\infty}=\sup_{0\leq x\leq 1}|f(x)|$ be the $l_{\infty}$-norm of $f$ defined as in the appendix, then $|\mu_t^N(f)|\leq \|f\|_{\infty}\frac{1}{N}\sum_{i=1}^NX_t^N(i)$. According to the definition of $\{X_t^N\}_{t\geq 0}$, $\sum_{i=1}^NX_t^N(i)\equiv \sum_{i=1}^NX_t^N(0)$ for any $t\geq 0$. Hence
\[
P\big(|\mu_t^N(f)|\geq M, \frac{1}{N}\sum_{i=1}^NX_0^N(i)\leq 2\int_0^1 \phi(x)dx\big)=0
\]
when $M>2\|f\|_{\infty}\int_0^1 \phi(x)dx$ and hence
\[
\lim_{M\rightarrow+\infty}\limsup_{N\rightarrow+\infty}P\big(|\mu_t^N(f)|\geq M\big)\leq \limsup_{N\rightarrow+\infty}P\big(\frac{1}{N}\sum_{i=1}^NX_0^N(i)>2\int_0^1 \phi(x)dx\big).
\]
Under Assumption (A), by Chebyshev's inequality,
\[
P\big(\frac{1}{N}\sum_{i=1}^NX_0^N(i)>2\int_0^1 \phi(x)dx\big)\leq \frac{1}{N\big(\int_0^1 \phi(x)dx+o(1)\big)^2} \rightarrow 0
\]
as $N\rightarrow 0$ and hence Condition (1) holds.

To check Condition (2), note that $|\mu_u^N(f)-\mu_{u-}^N(f)|\leq \frac{2\|f\|_{\infty}}{N}$ when $\{X_t^N\}_{t\geq 0}$ jumps at moment $u$,
hence, by strong Markov property, $|\mu_{\tau+s}^N(f)-\mu_\tau^N(f)|$ is stochastically dominated from above by $\frac{2\|f\|_{\infty}}{N}Y(2N\delta\|\lambda\|_{\infty}\int_0^1\phi(x)dx)$ conditioned on $\frac{1}{N}\sum_{i=1}^NX_0^N(i)\leq 2\int_0^1 \phi(x)dx$, where $\{Y(t)\}_{t\geq 0}$ is a Poisson process with rate one and $\|\lambda\|_{\infty}=\sup_{0\leq x,y\leq 1}\lambda(x,y)$. As a result, by Chebyshev's inequality,
\begin{align*}
&\lim_{\delta\rightarrow 0}\limsup_{N\rightarrow+\infty}\sup_{\tau\in \mathcal{T}, s\leq\delta}P\big(|\mu_{\tau+s}^N(f)-\mu_\tau^N(f)|>\epsilon\big)\\
&\leq \lim_{\delta\rightarrow 0}\limsup_{N\rightarrow+\infty}P\big(\frac{2\|f\|_{\infty}}{N}Y(2N\delta\|\lambda\|_{\infty}\int_0^1\phi(x)dx)>\epsilon\big)\\
&\text{\quad}+\limsup_{N\rightarrow+\infty}P\big(\frac{1}{N}\sum_{i=1}^NX_0^N(i)>2\int_0^1 \phi(x)dx\big)\\
&\leq \limsup_{\delta\rightarrow 0}\limsup_{N\rightarrow+\infty}e^{N[2\delta\|\lambda\|_{\infty}(e^{2\|f\|_{\infty}}-1)\int_0^1\phi(x)dx-\epsilon]}+\lim_{N\rightarrow+\infty}\frac{1}{N\big(\int_0^1 \phi(x)dx+o(1)\big)^2}\\
&=0
\end{align*}
and hence Condition (2) holds.

\qed

Now we give the proof of Theorem \ref{theorem 2.1 main result hydrolimit}.

\proof[Proof of Theorem \ref{theorem 2.1 main result hydrolimit}]
As we have shown in Lemma \ref{Existen and uniqueness}, Equation \eqref{equ 2.1 ODE on banach} has a unique weak solution.
Hence, by Lemma \ref{lemma 3.1}, we only need to show that any limit point $\widetilde{\mu}$ of subsequence of $\{\mu^N\}_{N\geq 1}$ satisfies Equation \eqref{equation 2.1 Mu with f} for any $f\in C[0,1]$.

For any $f\in C[0,1]$ and $N\geq 1$, we define
\[
M_t^N(f)=\mu_t^N(f)-\mu_0^N(f)-\int_0^t\mathcal{L}_N\big(\mu_s^N(f)\big)ds,
\]
then, by Dynkin's martingale formula, $\{M_t^N(f)\}_{t\geq 0}$ is a martingale with $\langle M^N(f)\rangle_t$ given by
\[
\langle M^N(f)\rangle_t=\int_0^t\Big(\mathcal{L}_N\left(\left(\mu_s^N(f)\right)^2\right)-2\mu_s^N(f)\mathcal{L}_N\left(\mu_s^N(f)\right)\Big)ds.
\]
According to the definition of $\mathcal{L}_N$,
\begin{align*}
&\mathcal{L}_N\left(\mu_s^N(f)\right)=\sum_{i=1}^N\sum_{j=1}^N\frac{X_s^N(i)\lambda(\frac{i}{N}, \frac{j}{N})}{N}\left(\mu_s^N(f)+\frac{f(\frac{j}{N})-f(\frac{i}{N})}{N}-\mu_s^N(f)\right) \\
&=-\frac{1}{N}\sum_{i=1}^NX_s^N(i)f(\frac{i}{N})\big[\frac{1}{N}\sum_{j=1}^N\lambda(\frac{i}{N}, \frac{j}{N})\big]+\frac{1}{N}\sum_{i=1}^NX_s^N(i)(\frac{1}{N}\sum_{j=1}^N\lambda(\frac{i}{N}, \frac{j}{N})f(\frac{j}{N}))\\
&=-\frac{1}{N}\sum_{i=1}^NX_s^N(i)f(\frac{i}{N})\big(P_2f(\frac{i}{N})+o(1)\big)+\frac{1}{N}\sum_{i=1}^NX_s^N(i)\big(P_1f(\frac{i}{N})+o(1)\big).
\end{align*}
Note that $o(1)$s in the above equation can be chosen uniformly according to the absolute continuity of $\lambda$.
According to the fact that $\sum_{i=1}^NX_t^N(i)$ are conserved and Chebyshev's inequality, under Assumption (A),
\[
\frac{1}{N}\sum_{i=1}^NX_s^N(i)=\frac{1}{N}\sum_{i=1}^NX_0^N(i)=\int_0^1\phi(y)dy+o_p(1).
\]
Hence, $\mathcal{L}_N\big(\mu_s^N(f)\big)=\mu_s^N\left((P_1-P_2)f\right)+o_p(1)$.
Similarly,
\[
\mu_0^N(f)=\int_0^1 f(y)\phi(y)dy+o_p(1)
\]
and
\begin{align*}
&\mathcal{L}_N\left(\left(\mu_s^N(f)\right)^2\right)-2\mu_s^N(f)\mathcal{L}_N\left(\mu_s^N(f)\right)=\sum_{i=1}^N\sum_{j=1}^N\frac{X_s^N(i)\lambda(\frac{i}{N}, \frac{j}{N})}{N^3}(f(\frac{j}{N})-f(\frac{i}{N}))^2\\
&\leq \frac{4\|f\|_{\infty}^2\|\lambda\|_{\infty}}{N^2}\sum_{i=1}^NX_s^N(i)=\frac{4\|f\|_{\infty}^2\|\lambda\|_{\infty}}{N}\left(\int_0^1\phi(y)dy+o_p(1)\right)=o_p(1).
\end{align*}
Then, $M_t^N(f)=o_p(1)$ by Doob's inequality. Consequently,
\begin{align*}
\mu_t^N(f)&=\mu_0^N(f)+\int_0^t\mathcal{L}_N\big(\mu_s^N(f)\big)ds+M_t^N(f)\\
&=\int_0^1f(y)\phi(y)dy+\int_0^t\mu_s^N\left((P_1-P_2)f\right)ds+o_p(1).
\end{align*}
Therefore, for any limit point $\widetilde{\mu}$ of subsequence of $\{\mu^N\}_{N\geq 1}$ and $f\in C[0,1]$,
\[
\widetilde{\mu}_t(f)=\int_0^1f(y)\phi(y)dy+\int_0^t\widetilde{\mu}_s\left((P_1-P_2)f\right)ds,
\]
i.e., $\widetilde{\mu}$ satisfies Equation \eqref{equation 2.1 Mu with f} and hence the proof is complete.

\qed

\section{Proof of Theorem \ref{theorem 2.2 main result CLT}}\label{section four}
In this section, we prove Theorem \ref{theorem 2.2 main result CLT}. First we show the tightness of $\{V^N\}_{N\geq 1}$.

\begin{lemma}\label{lemma 4.1}
Under Assumption (A), $\{V^N\}_{N\geq 1}$ is tight.
\end{lemma}

To prove Lemma \ref{lemma 4.1}, we need the following lemma.

\begin{lemma}\label{lemma 4.2}
There exists a constant $C_1<+\infty$ such that
\[
{\rm E}\left(\left(V_t^N(f)\right)^2\right)={\rm Cov}\left(V_t^N(f), V_t^N(f)\right)\leq C_1\|f\|_\infty^2
\]
\end{lemma}
for all $N\geq 1$ and any $f\in C[0, 1]$.

We first utilize Lemma \ref{lemma 4.2} to prove Lemma \ref{lemma 4.1}.

\proof[Proof of Lemma \ref{lemma 4.1}]

Similarly with the proof of Lemma \ref{lemma 3.1}, we only need to check that the following two conditions holds.

(1) $\lim_{M\rightarrow+\infty}\limsup_{N\rightarrow+\infty}P\big(|V_t^N(f)|\geq M\big)=0$ for any $t\geq 0$ and $f\in C[0,1]$.

(2) For any $\epsilon>0$ and $f\in C[0, 1]$,
\[
\lim_{\delta\rightarrow 0}\limsup_{N\rightarrow+\infty}\sup_{\tau\in \mathcal{T}, s\leq\delta}P\big(|V^N_{\tau+s}(f)-V^N_{\tau}(f)|>\epsilon\big)=0.
\]

By Dynkin's martingale formula, $V_t^N(f)=V_0^N(f)+\alpha_t^N+\beta_t^N$, where $\alpha_t^N=\int_0^t \big(\mathcal{L}_N+\partial_s\big)V_s^N(f)ds$ while $\{\beta_t^N\}_{t\geq 0}$ is a martingale
with
\begin{equation}\label{equ 4.1}
\langle\beta^N\rangle_t=\int_0^t\mathcal{L}_N\big((V_s^N(f))^2\big)-2V_s^N(f)\mathcal{L}_NV_s^N(f)ds.
\end{equation}
Then, to check conditions (1) and (2), we only need to check that the following four conditions hold.

(3) $\lim_{M\rightarrow+\infty}\limsup_{N\rightarrow+\infty}P\big(|\alpha_t^N|\geq M\big)=0$ for any $t\geq 0$ and $f\in C[0,1]$.

(4) For any $\epsilon>0$ and $f\in C[0, 1]$,
\[
\lim_{\delta\rightarrow 0}\limsup_{N\rightarrow+\infty}\sup_{\tau\in \mathcal{T}, s\leq\delta}P\big(|\alpha^N_{\tau+s}-\alpha^N_{\tau}|>\epsilon\big)=0.
\]

(5) $\lim_{M\rightarrow+\infty}\limsup_{N\rightarrow+\infty}P\big(|\beta_t^N|\geq M\big)=0$ for any $t\geq 0$ and $f\in C[0,1]$.

(6) For any $\epsilon>0$ and $f\in C[0, 1]$,
\[
\lim_{\delta\rightarrow 0}\limsup_{N\rightarrow+\infty}\sup_{\tau\in \mathcal{T}, s\leq\delta}P\big(|\beta^N_{\tau+s}-\beta^N_{\tau}|>\epsilon\big)=0.
\]

We first check Condition (3). According to the expression of the generator $\mathcal{L}_N$ of $\{X_t^N\}_{t\geq 0}$ given in Equation \eqref{equ 1.1 generator},
\[
\partial_t{\rm E}\left(X_t^N(i)\right)=-\frac{1}{N}\sum_{j=1}^N{\rm E}\left(X_t^N(i)\right)\lambda\left(\frac{i}{N},\frac{j}{N}\right)
+\frac{1}{N}\sum_{j=1}^N{\rm E}\left(X_t^N(j)\right)\lambda\left(\frac{j}{N}, \frac{i}{N}\right)
\]
for all $1\leq i\leq N$. As a result,
\begin{align}\label{equ alphat}
&\big(\mathcal{L}_N+\partial_s\big)V_s^N(f)\notag\\
&=\sum_{i=1}^N\sum_{j=1}^N\frac{X_s^N(i)\lambda(\frac{i}{N}, \frac{j}{N})}{N}\frac{f(\frac{j}{N})-f(\frac{i}{N})}{\sqrt{N}}-\frac{1}{\sqrt{N}}\sum_{i=1}^Nf\left(\frac{i}{N}\right)\partial_s{\rm E}\left(X_t^N(i)\right)\notag\\
&=V_s^N\big((P_1^N-P_2^N)f\big),
\end{align}
where
\[
(P_1^Nf)(x)=\frac{1}{N}\sum_{j=1}^N\lambda\left(x, \frac{j}{N}\right)f\left(\frac{j}{N}\right) \text{~and~}
(P_2^Nf)(x)=\frac{f(x)}{N}\sum_{j=1}^N\lambda\left(x, \frac{j}{N}\right)
\]
for $x\in [0, 1]$. Then, by Markov's inequality, H\"{o}lder inequality and Lemma \ref{lemma 4.2},
\begin{align*}
P\big(|\alpha_t^N|\geq M\big)&\leq \frac{1}{M^2}{\rm E}\left(\left(\int_0^t V_s^N\left((P_1^N-P_2^N)f\right)ds\right)^2\right)\\
&=\frac{1}{M^2}\int_0^t\int_0^t{\rm Cov}\left(V_s^N\left((P_1^N-P_2^N)f\right), V_u^N\left((P_1^N-P_2^N)f\right)\right)dsdu\\
&\leq\frac{1}{M^2}\int_0^t\int_0^t \sqrt{{\rm Cov}\left(V_s^N\left((P_1^N-P_2^N)f\right), V_s^N\left((P_1^N-P_2^N)f\right)\right)} \\
& \text{\quad\quad}\times\sqrt{{\rm Cov}\left(V_u^N\left((P_1^N-P_2^N)f\right), V_u^N\left((P_1^N-P_2^N)f\right)\right)} dsdu\\
&\leq \frac{1}{M^2}\int_0^t\int_0^t\sqrt{C_1\|(P_1^N-P_2^N)f\|_\infty^2}\sqrt{C_1\|(P_1^N-P_2^N)f\|_\infty^2}dsdu\\
&=\frac{C_1t^2}{M^2}\|(P_1^N-P_2^N)f\|_\infty^2\leq \frac{4C_1t^2}{M^2}\|\lambda\|_{\infty}^2\|f\|^2_{\infty}
\end{align*}
and hence Condition (3) holds.

Then we check Condition (4). According to Equation \eqref{equ alphat},
\[
\alpha^N_{\tau+s}-\alpha^N_{\tau}=\int_\tau^{\tau+s} V_u^N\big((P_1^N-P_2^N)f\big) du.
\]
Therefore, by H\"{o}lder inequality, Markov's inequality and Lemma \ref{lemma 4.2}, for $s\leq \delta$ and $\tau\in \mathcal{T}$,
\begin{align*}
&P\big(|\alpha^N_{\tau+s}-\alpha^N_{\tau}|>\epsilon\big)\\
& \leq \frac{1}{\epsilon}{\rm E}\left(\int_\tau^{\tau+s} |V_u^N\big((P_1^N-P_2^N)f\big)| du\right)\\
&=\frac{1}{\epsilon}{\rm E}\left(\int_0^{T+\delta} |V_u^N\big((P_1^N-P_2^N)f\big)|1_{\{\tau\leq u\leq \tau+s\}}du\right)\\
&=\frac{1}{\epsilon}\int_0^{T+\delta}{\rm E}\left(  |V_u^N\big((P_1^N-P_2^N)f\big)|1_{\{\tau\leq u\leq \tau+s\}}\right)du \\
&\leq\frac{1}{\epsilon}\int_0^{T+\delta}\sqrt{{\rm Cov}\left(V_u^N\big((P_1^N-P_2^N)f\big), V_u^N\big((P_1^N-P_2^N)f\big)\right)}\sqrt{P(\tau\leq u\leq \tau+s)}du\\
&\leq \frac{\sqrt{C_1}\|(P_1^N-P_2^N)f\|_\infty}{\epsilon}\int_0^{T+\delta}\sqrt{P(\tau\leq u\leq \tau+s)} du\\
&\leq \frac{2\sqrt{C_1(T+\delta)}\|\lambda\|_{\infty}\|f\|_\infty}{\epsilon}\sqrt{\int_0^{T+\delta}P(\tau\leq u\leq \tau+s)du}\\
&=\frac{2\sqrt{C_1(T+\delta)}\|\lambda\|_{\infty}\|f\|_\infty}{\epsilon}\sqrt{{\rm E}\left(\int_0^{T+\delta}1_{\{\tau\leq u\leq \tau+s\}}du\right)}\\
&=\frac{2\sqrt{C_1(T+\delta)s}\|\lambda\|_{\infty}\|f\|_\infty}{\epsilon}\leq\frac{2\sqrt{C_1(T+\delta)\delta}\|\lambda\|_\infty\|f\|_\infty}{\epsilon}
\end{align*}
and hence Condition (4) holds.

Then we check Condition (5). According to the definition of $\mathcal{L}_N$,
\begin{align*}
&\mathcal{L}_N\big((V_s^N(f))^2\big)-2V_s^N(f)\mathcal{L}_NV_s^N(f)=\frac{1}{N^2}\sum_{i=1}^N\sum_{j=1}^NX_s^N(i)\lambda(\frac{j}{N}, \frac{i}{N})\big(f(\frac{j}{N})-f(\frac{i}{N})\big)^2 \\
&\leq \left(\frac{4}{N}\sum_{i=1}^NX_s^N(i)\right)\|\lambda\|_{\infty}\|f\|_{\infty}^2=\left(\frac{4}{N}\sum_{i=1}^NX_0^N(i)\right)\|\lambda\|_{\infty}\|f\|_{\infty}^2.
\end{align*}
Hence,
\begin{align*}
{\rm E}\left((\beta_t^N)^2\right)={\rm E}\big(\langle \beta^N\rangle_t\big)&\leq \frac{4t}{N}\left(\sum_{i=1}^N\phi\left(\frac{i}{N}\right)\right)\|\lambda\|_{\infty}\|f\|_{\infty}^2\\
&=\left(4t\int_0^1\phi(x)dx+o(1)\right)\|\lambda\|_{\infty}\|f\|_{\infty}^2
\end{align*}
and Condition (5) follows from Chebyshev's inequality.

 At last we check Condition (6). Since ${\rm E}\big[|\beta_{\tau+s}^N-\beta_\tau^N|^2\big]={\rm E}\big(\langle \beta^N\rangle_{\tau+s}-\langle \beta^N\rangle_\tau\big)$ while
\begin{align*}
&\langle \beta^N\rangle_{\tau+s}-\langle \beta^N\rangle_\tau\leq \int_{\tau}^{\tau+s} \left(\frac{4}{N}\sum_{i=1}^NX_0^N(i)\right)\|\lambda\|_{\infty}\|f\|_{\infty}^2du\\
&\leq \left(\frac{4\delta}{N}\sum_{i=1}^NX_0^N(i)\right)\|\lambda\|_{\infty}\|f\|_{\infty}^2
\end{align*}
for any $s\leq \delta$,
\[
{\rm E}\big[|\beta_{\tau+s}^N-\beta_\tau^N|^2\big]
\leq 4\delta\left(\int_0^1\phi(x)dx+o(1)\right)\|\lambda\|_{\infty}\|f\|_{\infty}^2.
\]
As a result, Condition (6) follows from Chebyshev's inequality.

Since Conditions (3)-(6) all hold, the proof is complete.

\qed

Now we give proofs of Lemma \ref{lemma 4.2}.

\proof[Proof of Lemma \ref{lemma 4.2}]
By direct calculation,
\[
{\rm Cov}\left(V_t^N(f), V_t^N(f)\right)\leq {\rm \uppercase\expandafter{\romannumeral1}}+{\rm \uppercase\expandafter{\romannumeral2}},
\]
where
\[
{\rm \uppercase\expandafter{\romannumeral1}}=\frac{\|f\|_{\infty}^2}{N}\sum_{i=1}^N{\rm E}\left(\left(X_t^N(i)\right)^2\right)
\]
and
\[
{\rm \uppercase\expandafter{\romannumeral2}}=\frac{1}{N}\sum_{i=1}^N\sum_{j\neq i}^Nf\left(\frac{i}{N}\right)f\left(\frac{j}{N}\right)
{\rm Cov}\left(X_t^N(i), X_t^N(j)\right).
\]
Conditioned on $X_0^N$, since each ball jumps to the $i$th box at rate at most $\frac{\|\lambda\|_{\infty}}{N}$, $X_t^N(i)-X_0^N(i)$ is stochastically dominated from above by a Poisson process $\{K_t\}_{t\geq 0}$ with rate
\[
\|\lambda\|_\infty\frac{1}{N}\sum_{j=1}^NX_0^N(j).
\]
Hence,
\begin{align*}
&{\rm E}\left(\left(X_t^N(i)\right)^2\big|X_0^N\right) \leq \left(X_0^N(i)\right)^2+2X_0^N(i){\rm E}K_t+{\rm E}\left(K_t^2\right)\\
&=\left(X_0^N(i)\right)^2+2X_0^N(i)\|\lambda\|_{\infty}\frac{t}{N}\sum_{j=1}^NX_0^N(j)+\|\lambda\|_{\infty}\frac{t}{N}\sum_{j=1}^NX_0^N(j)\\
&\text{\quad\quad}+\frac{t^2\|\lambda\|_{\infty}^2}{N^2}\left(\sum_{j=1}^NX_0^N(j)\right)^2
\end{align*}
and hence
\begin{align}\label{equ 4.2 upper bound of Roma1}
&{\rm E}\left(\left(X_t^N(i)\right)^2\right)\\
&\leq \|\phi\|_{\infty}^2\left(1+2\|\lambda\|_\infty T+T^2\|\lambda\|^2_{\infty}\right)
+\|\phi\|_{\infty}\left(1+3\|\lambda\|_{\infty}T+T^2\|\lambda\|^2_{\infty}\right) \notag
\end{align}
for all $0\leq t\leq T$ according to Assumption (A).

For $i\neq j$, ${\rm Cov}\left(X_t^N(i), X_t^N(j)\right)={\rm \uppercase\expandafter{\romannumeral3}}+{\rm \uppercase\expandafter{\romannumeral4}}$, where
\[
{\rm \uppercase\expandafter{\romannumeral3}}={\rm E}\left({\rm E}\left(X_t^N(i)X_t^N(j)\big|X_0^N\right)-{\rm E}\left(X_t^N(i)\big|X_0^N\right){\rm E}\left(X_t^N(j)\big|X_0^N\right)\right)
\]
and
\[
{\rm \uppercase\expandafter{\romannumeral4}}={\rm Cov}\left({\rm E}\left(X_t^N(i)\big|X_0^N\right), {\rm E}\left(X_t^N(j)\big|X_0^N\right)\right).
\]
For a given ball initially in box $i$, we use $p_{ij}(t)$ to denote the probability that this ball is in box $j$ at moment $t$. Then,
\[
{\rm E}\left(X_t^N(i)\big|X_0^N\right)=\sum_{k=1}^NX^N_0(k)p_{ki}(t).
\]
Hence, by Assumption (A),
\begin{align*}
{\rm \uppercase\expandafter{\romannumeral4}}&=\sum_{k=1}^N\sum_{l=1}^N{\rm Cov}\left(X_0^N(k), X_0^N(l)\right)p_{ki}(t)p_{lj}(t)\\
&=\sum_{k=1}^N{\rm Cov}\left(X_0^N(k), X_0^N(k)\right)p_{ki}(t)p_{kj}(t)=\sum_{k=1}^N\phi\left(\frac{k}{N}\right)p_{ki}(t)p_{kj}(t).
\end{align*}
For $1\leq i,j\leq N$ and $1\leq l\leq X_0^N(i)$, we use $A_t^i(l, j)$ to denote the event that the $l$-th ball initially in box $i$ is in box $j$ at moment $t$. Then, $X_t^N(i)=\sum_{k=1}^N\sum_{l=1}^{X_0^N(k)}1_{A_t^k(l, i)}$ and
\begin{align*}
&{\rm E}\left(X_t^N(i)X_t^N(j)\big|X_0^N\right)-{\rm E}\left(X_t^N(i)\big|X_0^N\right){\rm E}\left(X_t^N(j)\big|X_0^N\right)\\
&=\sum_{k=1}^N\sum_{l=1}^{X_0^N(k)}\sum_{(u,s)\neq (k,l)}P\left(A^k_t(l,i)\right)P\left(A^u_t(s, j)\right)\\
&\text{\quad\quad}-\sum_{k=1}^N\sum_{l=1}^{X_0^N(k)}\sum_{u=1}^N\sum_{s=1}^{X_0^N(u)}P\left(A^k_t(l,i)\right)P\left(A^u_t(s, j)\right)\\
&=-\sum_{k=1}^N\sum_{l=1}^{X_0^N(k)}P\left(A^k_t(l,i)\right)P\left(A^k_t(l, j)\right)=-\sum_{k=1}^NX_0^N(k)p_{ki}(t)p_{kj}(t).
\end{align*}
Therefore, by Assumption (A),
\[
{\rm \uppercase\expandafter{\romannumeral3}}=-\sum_{k=1}^N\phi\left(\frac{k}{N}\right)p_{ki}(t)p_{kj}(t)=-{\rm \uppercase\expandafter{\romannumeral4}}
\]
and hence ${\rm \uppercase\expandafter{\romannumeral2}}=0$. As a result, Lemma \ref{lemma 4.2} follows from Equation \eqref{equ 4.2 upper bound of Roma1}.

\qed

At last, we give the proof of Theorem \ref{theorem 2.2 main result CLT}.

\proof[Proof of Theorem \ref{theorem 2.2 main result CLT}]

By Lemma \ref{lemma 4.1}, we only need to show that any limit point $\widetilde{V}$ of subsequence of $\{V^N\}_{N\geq 1}$ is a solution of Equation \eqref{equ 2.1}. Equivalently, as we have introduced in Section \ref{section two}, we need to check that $\widetilde{V}$ has those three properties corresponding to the definition of the unique solution of Equation \eqref{equ 2.1}. The first property that $\widetilde{V}(H)$ is continuous for any $H\in C[0, 1]$ follows from the fact that
\[
\sup_{0\leq t\leq T}|V_t^N(H)-V_{t-}^N(H)|\leq \frac{2\|H\|_{\infty}}{\sqrt{N}}.
\]
The third property that $\widetilde{V}_0(H)$ follows $\mathbb{N}\left(0, \int_0^1H^2(x)\phi(x) dx\right)$ is a consequence of Assumption (A). Hence, we only need to check that $\widetilde{V}_0(H)$ has the second property, i.e.,
\begin{align*}
\Bigg\{G\left(\widetilde{V}_t(H)\right)-G\left(\widetilde{V}_0(H)\right)&-\int_0^tG^\prime\left(\widetilde{V}_s(H)\right)\widetilde{V}_s\left((P_1-P_2)H\right)ds \\
&-\frac{1}{2}\int_0^tG^{\prime\prime}\left(\widetilde{V}_s(H)\right)\|b_sH\|_2^2ds\Bigg\}_{0\leq t\leq T}
\end{align*}
is a martingale for $H\in C[0,1]$ and $G\in C_c^{\infty}(\mathbb{R})$.

For any $G\in C_c^{\infty}(\mathbb{R})$ and $H\in C[0,1]$, we define
\[
\mathcal{M}_t^N(G,H)=G\big(V_t^N(H)\big)-G\big(V_0^N(H)\big)-\int_0^t\big(\mathcal{L}_N+\partial_s\big)G\big(V_s^N(H)\big)ds,
\]
then, by Dynkin's martingale formula, $\{\mathcal{M}_t^N(G, H)\}_{t\geq 0}$ is a martingale. According to the definition of $\mathcal{L}_N$ and Taylor's expansion up to the third order,
\begin{align}\label{equ 4.4 generatorandpartial}
&\big(\mathcal{L}_N+\partial_s\big)G\big(V_s^N(H)\big)\notag \\
&=\sum_{i=1}^N\sum_{j=1}^N\frac{X_s^N(i)\lambda\big(\frac{i}{N}, \frac{j}{N}\big)}{N}\Big[G\left(V_s^N(H)
+\frac{H(\frac{j}{N})-H(\frac{i}{N})}{\sqrt{N}}\right)-G\big(V_s^N(H)\big)\Big]\notag\\
&\text{\quad}-\frac{G^\prime\big(V_s^N(H)\big)}{\sqrt{N}}\sum_{i=1}^N\frac{H\left(\frac{i}{N}\right)}{N}\big[\sum_{j=1}^N{\rm E}\left(X_s^N(j)\right)\lambda(\frac{j}{N},\frac{i}{N})-\sum_{j=1}^N{\rm E}\left(X_s^N(i)\right)\lambda(\frac{i}{N},\frac{j}{N})\big]\notag\\
&=\sum_{i=1}^N\sum_{j=1}^N\frac{X_s^N(i)\lambda(\frac{i}{N}, \frac{j}{N})}{N}G^\prime\big(V_s^N(H)\big)\frac{H(\frac{j}{N})-H(\frac{i}{N})}{\sqrt{N}}\notag\\
&\text{\quad}+\frac{1}{2}\sum_{i=1}^N\sum_{j=1}^N\frac{X_s^N(i)\lambda(\frac{i}{N}, \frac{j}{N})}{N}G^{\prime\prime}\big(V_s^N(H)\big)\frac{\big(H(\frac{j}{N})-H(\frac{i}{N})\big)^2}{N}
+\xi_s^N\notag\\
&\text{\quad}-\frac{G^\prime\big(V_s^N(H)\big)}{\sqrt{N}}\sum_{i=1}^N\frac{H\left(\frac{i}{N}\right)}{N}\big[\sum_{j=1}^N{\rm E}\left(X_s^N(j)\right)\lambda(\frac{j}{N},\frac{i}{N})-\sum_{j=1}^N{\rm E}\left(X_s^N(i)\right)\lambda(\frac{i}{N},\frac{j}{N})\big]\notag\\
&=G^\prime\big(V_s^N(H)\big)V_s^N\left(\left(P_1^N-P_2^N\right)H\right)+\xi_s^N\\
&\text{\quad}+\frac{1}{2}\sum_{i=1}^N\sum_{j=1}^N\frac{X_s^N(i)\lambda(\frac{i}{N}, \frac{j}{N})}{N}G^{\prime\prime}\big(V_s^N(H)\big)\frac{\big(H(\frac{j}{N})-H(\frac{i}{N})\big)^2}{N},\notag
\end{align}
where $\xi_s^N$ is bounded from above by
\begin{align*}
&\frac{\|G^{\prime\prime\prime}\|_{\infty}}{6}\sum_{i=1}^N\sum_{j=1}^N\frac{X_s^N(i)\lambda(\frac{i}{N}, \frac{j}{N})}{N}\frac{\left|H(\frac{j}{N})-H(\frac{i}{N})\right|^3}{N^{\frac{3}{2}}}\\
&\leq \frac{4\|G^{\prime\prime\prime}\|_{\infty}\|\lambda\|_{\infty}\|H\|_{\infty}}{3N^{\frac{1}{2}}}\left(\frac{1}{N}\sum_{i=1}^NX_0^N(i)\right)\\
&=\frac{4\|G^{\prime\prime\prime}\|_{\infty}\|\lambda\|_{\infty}\|H\|_{\infty}}{3N^{\frac{1}{2}}}\left(\int_0^1\phi(x)dx+o_p(1)\right)=o_p(1).
\end{align*}
According to the definition of $P_1^N, P_2^N$ and the fact that $\lambda\in C^{1,1}\left([0,1]\times[0, 1]\right)$, it is easy to check that $\|(P_1^N-P_2^N)H-(P_1-P_2)H\|_{\infty}=\frac{O(1)}{N}$ while
\[
V_s^N\left(\left(P_1^N-P_2^N\right)H\right)-V_s^N\left(\left(P_1-P_2\right)H\right)=o_p(1)
\]
and this $o_p(1)$ can be chosen uniformly for $0\leq s\leq T$. Since $\frac{1}{N}\sum_{i=1}^NX_s^N(i)=\frac{1}{N}\sum_{i=1}^NX_0^N(i)=\int_0^1\phi(x)dx+o_p(1)$
and
\begin{align*}
&\sup_{1\leq i\leq N}\left|\frac{\sum_{j=1}^N\lambda(\frac{i}{N},\frac{j}{N})\left(H(\frac{j}{N})-H(\frac{i}{N})\right)^2}{N}-\int_0^1\lambda(\frac{i}{N}, y)\left(H(y)-H\left(\frac{i}{N}\right)\right)^2dy\right|\\
&=o(1),
\end{align*}
by Theorem \ref{theorem 2.1 main result hydrolimit},
\begin{align*}
&\sum_{i=1}^N\sum_{j=1}^N\frac{X_s^N(i)\lambda(\frac{i}{N}, \frac{j}{N})}{N}\frac{\big(H(\frac{j}{N})-H(\frac{i}{N})\big)^2}{N}\\
&=\int_0^1\rho(s,x)\left[\int_0^1\lambda(x,y)\left(H(y)-H(x)\right)^2dy\right]dx+o_p(1)=\|b_sH\|_2^2+o_p(1)
\end{align*}
and this $o_p(1)$ can be chosen uniformly for $0\leq s\leq T$.

Consequently,
\[
\big(\mathcal{L}_N+\partial_s\big)G\big(V_s^N(H)\big)=G^\prime\big(V_s^N(H)\big)V_s^N\big((P_1-P_2)H\big)+\frac{1}{2}G^{\prime\prime}\big(V_s^N(H)\big)\|b_sH\|_2^2+o_p(1)
\]
and
\begin{align}\label{equ 4.3 Mt}
&\mathcal{M}_t^N(G, H)=G\big(V_t^N(H)\big)-G\big(V_0^N(H)\big)-\int_0^t G^\prime\big(V_s^N(H)\big)V_s^N\big((P_1-P_2)H\big)ds\notag\\
&-\int_0^t\frac{1}{2}G^{\prime\prime}\big(V_s^N(H)\big)\|b_sH\|_2^2ds+o_p(1).
\end{align}
By Equation \eqref{equ 4.3 Mt}, a subsequence of $\{\mathcal{M}_t^N(G, H)\}_{0\leq t\leq T}$ converges weakly to
\begin{align*}
\Bigg\{G\left(\widetilde{V}_t(H)\right)-G\left(\widetilde{V}_0(H)\right)&-\int_0^tG^\prime\left(\widetilde{V}_s(H)\right)\widetilde{V}_s\left((P_1-P_2)H\right)ds \\
&-\frac{1}{2}\int_0^tG^{\prime\prime}\left(\widetilde{V}_s(H)\right)\|b_sH\|_2^2ds\Bigg\}_{0\leq t\leq T}.
\end{align*}
Since $\{\mathcal{M}_t^N(G, H)\}_{t\geq 0}$ is a martingale for each $N\geq 1$, according to Theorem 5.3 of \cite{Whitt2007}, we only need to show that $\{\mathcal{M}_t^N(G, H)\}_{N\geq 1}$ are uniformly integrable for each $t\geq 0$ to complete this proof. Since $G\in C_c^{\infty}(\mathbb{R})$, we only need to show that
\begin{equation}\label{equ 4.5}
\sup_{0\leq s\leq T, N\geq 1}{\rm E}\left(\left(\big(\mathcal{L}_N+\partial_s\big)G\big(V_s^N(H)\big) \right)^2\right)<+\infty
\end{equation}
to derive the uniform integrability of $\{\mathcal{M}_t^N(G, H)\}_{N\geq 1}$. According to Equation \eqref{equ 4.4 generatorandpartial}, Lemma \ref{lemma 4.2} and the upper bound of $\xi_s^N$ given above, Equation \eqref{equ 4.5} follows from the fact that
\[
\sup_{N\geq 1}{\rm E}\left(\left(\frac{1}{N}\sum_{i=1}^NX_0^N(i)\right)^2\right)\leq \|\phi\|_\infty^2+\|\phi\|_{\infty}
\]
given by Assumption (A) and the proof is complete.

\qed

\section{Proof of Equation \eqref{equ LDP upper}}\label{section five}

In this section, we give the proof of Equation \eqref{equ LDP upper}. As we have mentioned in Section \ref{section two}, an exponential martingale is crucial for our proof, so we first introduce this martingale. For any $H\in C^{1, 0}\left([0, T]\times X\right)$ and $0\leq t\leq T$, we use $\mathcal{U}_t^N(H)$ to denote
\[
H(t, X_t^N)-H(0, X_0^N)-\int_0^t\left(\partial_s+\mathcal{L}_N\right)H(s, X_s^N)ds.
\]
Then, by Dynkin's martingale formula, $\{\mathcal{U}_t^N(H)\}_{0\leq t\leq T}$ is a martingale. For any $G\in C^{1,0}\left([0, T]\times [0, 1]\right)$ and $x\in X$, we define
\[
f^N_G(t, x)=e^{\sum_{i=1}^Nx(i)G_t(\frac{i}{N})}
\]
and hence
\[
f^N_G(t, X_t^N)=e^{N\mu_t^N(G_t)}.
\]
For any $0\leq t\leq T$ and $G\in C^{1,0}\left([0, T]\times [0, 1]\right)$, we define
\[
\Lambda_t^N(G)=\frac{f_G^N(t, X_t^N)}{f_G^N(0, X_0^N)}e^{-\int_0^t{\frac{(\partial_s+\mathcal{L}_N)f_G^N(s, X_s^N)}{f_G^N(s, X_s^N)}}ds}.
\]
By the definition of $\mathcal{L}_N$ and direct calculation,
\begin{equation}\label{equ 5.1 Lambda}
\Lambda_t^N(G)=\exp\left\{N\mu_t^N(G_t)-N\mu_0^N(G_0)-N\int_0^t\mu_s^N\left((\partial_s+\mathcal{B}^N)G_s\right)ds\right\},
\end{equation}
where
\[
(\mathcal{B}^Ng)(x)=\frac{1}{N}\sum_{j=1}^N\lambda(x, \frac{j}{N})\left(e^{g(\frac{j}{N})-g(x)}-1\right)
\]
for any $g\in C[0, 1]$ and $x\in [0, 1]$. We have the following lemma.

\begin{lemma}\label{lemma LambdaisMartingale}
For any $G\in C^{1,0}\left([0, T]\times [0, 1]\right)$, $\{\Lambda_t^N(G)\}_{0\leq t\leq T}$ is a martingale.
\end{lemma}

\proof

Let $V_t^N(G)=\frac{1}{f_G^N(0, X_0^N)}e^{-\int_0^t{\frac{(\partial_s+\mathcal{L}_N)f_G^N(s, X_s^N)}{f_G^N(s, X_s^N)}}ds}$, then, by It\^{o}'s formula,
\begin{equation}\label{equ 5.2 dLambda}
d\Lambda_t^N(G)=V_t^N(G)d\mathcal{U}_t^N(f_G^N).
\end{equation}
As a result, $\{\Lambda_t^N(G)\}_{0\leq t\leq T}$ is a local martingale. By Equation \eqref{equ 5.1 Lambda},
\begin{equation}\label{equ 5.2 two}
\sup_{0\leq t\leq T}\Lambda_t^N(G)\leq \exp\left\{\sum_{i=1}^NX_0^N(i)\left(2+\|\lambda\|_\infty T\left(e^{2\|g\|_\infty}-1\right)+\|\partial_tG_t\|_\infty\right)\right\}.
\end{equation}
As a result, $\{\Lambda_t^N(G)\}_{0\leq t\leq T}$ being a martingale follows from the fact that
\[
{\rm E}\left(e^{\theta\sum_{i=1}^NX_0^N(i)}\right)<+\infty
\]
for any $\theta>0$ given by Assumption (A).

\qed

We will first show that Equation \eqref{equ LDP upper} holds for any compact $C$. A crucial step is to bound the error from above when we replace $\mu_t^N(\mathcal{B}^NG_t)$ by $\mu_t^N(\mathcal{B}G_t)$. So we need the following lemma.

\begin{lemma}\label{lemma 5.1 replacement}
For any $G\in C^{1, 0}\left([0, T]\times [0, 1]\right)$ and $\epsilon>0$,
\[
\limsup_{N\rightarrow+\infty}\frac{1}{N}\log P\left(\sup_{0\leq t\leq T}\left|\mu_t^N(\mathcal{B}G_t)-\mu_t^N(\mathcal{B}^NG_t)\right|>\epsilon\right)=-\infty.
\]
\end{lemma}

\proof

According to the uniform continuities of $\lambda$ and $G$, for any $\delta>0$, there exists $N(\delta)<+\infty$ such that
\[
\sup_{0\leq t\leq T, 0\leq x\leq 1}\left|\mathcal{B}G_t(x)-\mathcal{B}^NG_t(x)\right|<\delta
\]
when $N\geq N(\delta)$. Hence, when $N\geq N(\delta)$,
\[
\sup_{0\leq t\leq T}\left|\mu_t^N(\mathcal{B}^NG_t)-\mu_t^N(\mathcal{B}G_t)\right|\leq \frac{\delta}{N}\sum_{i=1}^NX_0^N(i).
\]
According to Assumption (A) and Chebyshev's inequality,
\begin{align}\label{equ 5.3}
P\left(\frac{\delta}{N}\sum_{i=1}^NX_0^N(i)\geq \epsilon\right)&\leq \exp\left\{-\frac{N\epsilon}{\delta}+(e-1)\sum_{i=1}^N\phi(\frac{i}{N})\right\}\notag\\
&=\exp\left\{N\left(-\frac{\epsilon}{\delta}+(e-1)\int_0^1\phi(x)dx+o(1)\right)\right\}.
\end{align}
Since $\delta$ is arbitrary, Lemma \ref{lemma 5.1 replacement} follows from Equation \eqref{equ 5.3}.

\qed

According to Lemmas \ref{lemma LambdaisMartingale} and \ref{lemma 5.1 replacement}, we have the following conclusion.

\begin{lemma}\label{lemma 5.3 UBholds for compact}
Under Assumptions (A), for any compact set $C\subseteq D\big([0, T], \mathcal{S}\big)$,
\begin{equation*}
\limsup_{N\rightarrow+\infty}\frac{1}{N}\log P\left(\mu^N\in C\right)\leq -\inf_{\mu\in C}\left(I_{ini}(\mu_0)+I_{dyn}(\mu)\right).
\end{equation*}
\end{lemma}

\proof

For any $\epsilon>0$, by Lemma \ref{lemma 5.1 replacement},
\begin{align}\label{equ 5.5}
&\limsup_{N\rightarrow+\infty}\frac{1}{N}\log P\left(\mu^N\in C\right)\\
&=\limsup_{N\rightarrow+\infty}\frac{1}{N}\log P\left(\mu^N\in C, \sup_{0\leq t\leq T}\left|\mu_t^N(\mathcal{B}G_t)-\mu_t^N(\mathcal{B}^NG_t)\right|\leq \epsilon\right).\notag
\end{align}
By Lemma \ref{lemma LambdaisMartingale} and Assumption (A), for any $H\in C[0, 1]$ and $G\in C^{1, 0}\left([0, T]\times [0, 1]\right)$,
\[
{\rm E}\left(e^{N\mu_0^N(H)}\Lambda_T^N(G)\right)={\rm E}\left(e^{N\mu_0^N(H)}\right)=e^{N\left(\int_0^1\phi(x)\left(e^{H(x)}-1\right)dx+o(1)\right)}.
\]
For simplicity, we use $r^N_{_{C, \epsilon}}$ to denote the event
\[
\left\{\mu^N\in C, \sup_{0\leq t\leq T}\left|\mu_t^N(\mathcal{B}G_t)-\mu_t^N(\mathcal{B}^NG_t)\right|\leq \epsilon\right\}.
\]
On $r^N_{_{C, \epsilon}}$,
\begin{align*}
&e^{N\mu_0^N(H)}\Lambda_T^N(G)\\
&\geq\exp\left\{N\mu_0^N(H)+N\mu_T^N(G_T)-N\mu_0^N(G_0)-N\int_0^T\mu_s^N\left((\partial_s+\mathcal{B})G_s\right)ds-NT\epsilon\right\}\\
&\geq\exp\left\{N\inf_{\mu\in C}\left\{\mu_0(H)+\mu_T(G_T)-\mu_0(G_0)-\int_0^T\mu_s\left((\partial_s+\mathcal{B})G_s\right)ds\right\}-NT\epsilon\right\}.
\end{align*}
Therefore,
\begin{align*}
&{\rm E}\left(e^{N\mu_0^N(H)}\right)=e^{N\left(\int_0^1\phi(x)\left(e^{H(x)}-1\right)dx+o(1)\right)}\\
&={\rm E}\left(e^{N\mu_0^N(H)}\Lambda_T^N(G)\right)\geq {\rm E}\left(e^{N\mu_0^N(H)}\Lambda_T^N(G)1_{r^N_{_{C, \epsilon}}}\right)\\
&\geq P(r^N_{_{C, \epsilon}})\\
&\times \exp\left\{N\inf_{\mu\in C}\left\{\mu_0(H)+\mu_T(G_T)-\mu_0(G_0)-\int_0^T\mu_s\left((\partial_s+\mathcal{B})G_s\right)ds\right\}-NT\epsilon\right\}
\end{align*}
and hence
\begin{align*}
&\limsup_{N\rightarrow+\infty}\frac{1}{N}\log P\left(r^N_{_{C, \epsilon}}\right)\leq \epsilon T-\inf_{\mu\in C}\Bigg\{\mu_0(H)-\int_0^1\phi(x)\left(e^{H(x)}-1\right)dx\\
&\text{\quad}+\mu_T(G_T)-\mu_0(G_0)-\int_0^T\mu_s\left((\partial_s+\mathcal{B})G_s\right)ds\Bigg\}.
\end{align*}
By Equation \eqref{equ 5.5},
\begin{align*}
&\limsup_{N\rightarrow+\infty}\frac{1}{N}\log P\left(\mu^N\in C\right)\leq \epsilon T-\inf_{\mu\in C}\Bigg\{\mu_0(H)-\int_0^1\phi(x)\left(e^{H(x)}-1\right)dx\\
&\text{\quad}+\mu_T(G_T)-\mu_0(G_0)-\int_0^T\mu_s\left((\partial_s+\mathcal{B})G_s\right)ds\Bigg\}.
\end{align*}
Since $\epsilon, G, H$ are arbitrary,
\begin{align*}
&\limsup_{N\rightarrow+\infty}\frac{1}{N}\log P\left(\mu^N\in C\right)\leq -\sup_{H\in C[0, 1],\atop G\in C^{1, 0}\left([0, T]\times [0, 1]\right)}\inf_{\mu\in C}\Bigg\{\mu_0(H)\\
&-\int_0^1\phi(x)\left(e^{H(x)}-1\right)dx
+\mu_T(G_T)-\mu_0(G_0)-\int_0^T\mu_s\left((\partial_s+\mathcal{B})G_s\right)ds\Bigg\}.
\end{align*}
Since $C$ is compact, according to the minimax theorem given in \cite{Sion1958},
\begin{align*}
&\sup_{H\in C[0, 1],\atop G\in C^{1, 0}\left([0, T]\times [0, 1]\right)}\inf_{\mu\in C}\Bigg\{\mu_0(H)-\int_0^1\phi(x)\left(e^{H(x)}-1\right)dx\\
&\text{\quad}+\mu_T(G_T)-\mu_0(G_0)-\int_0^T\mu_s\left((\partial_s+\mathcal{B})G_s\right)ds\Bigg\}\\
&=\inf_{\mu\in C}\sup_{H\in C[0, 1],\atop G\in C^{1, 0}\left([0, T]\times [0, 1]\right)}\Bigg\{\mu_0(H)-\int_0^1\phi(x)\left(e^{H(x)}-1\right)dx\\
&\text{\quad}+\mu_T(G_T)-\mu_0(G_0)-\int_0^T\mu_s\left((\partial_s+\mathcal{B})G_s\right)ds\Bigg\}\\
&=\inf_{\mu\in C}\Bigg\{\sup_{H\in C[0, 1]}\left\{\mu_0(H)-\int_0^1\phi(x)\left(e^{H(x)}-1\right)dx\right\}\\
&\text{\quad}+\sup_{G\in C^{1, 0}\left([0, T]\times [0, 1]\right)}\left\{\mu_T(G_T)-\mu_0(G_0)-\int_0^T\mu_s\left((\partial_s+\mathcal{B})G_s\right)ds\right\}\Bigg\}\\
&=\inf_{\mu\in C}\left(I_{ini}(\mu_0)+I_{dyn}(\mu)\right)
\end{align*}
and the proof is complete.

\qed

At last we give the proof of Equation \eqref{equ LDP upper} for all closed $C$.

\proof[Proof of Equation \eqref{equ LDP upper}]

By Lemma \ref{lemma 5.3 UBholds for compact}, we only need to show that $\{\mu^N\}_{N\geq 1}$ are exponential tight, i.e., for any $\epsilon>0$, there exists compact $K=K(\epsilon)\subset D\big([0, T], \mathcal{S}\big)$ such that
\[
\sup_{N\geq 1}\left\{\left(P(\mu^N\not\in K)\right)^{\frac{1}{N}}\right\}<\epsilon.
\]
According to the criteria given in \cite{Puhalskii1994}, we only need to show that
\begin{equation}\label{equ 5.5 EtightC1}
\limsup_{M\rightarrow+\infty}\limsup_{N\rightarrow+\infty}\frac{1}{N}\log P\left(\sup_{0\leq t\leq T}\left|\mu_t^N(H)\right|>M\right)=-\infty
\end{equation}
for any $H\in C[0, 1]$ and
\begin{equation}\label{equ 5.6 EtightC2}
\limsup_{\delta\rightarrow0}\limsup_{N\rightarrow+\infty}\frac{1}{N}\log\sup_{\tau\in \mathcal{T}}P\left(\sup_{0\leq t\leq \delta}\left|\mu_{t+\tau}^N(H)-\mu_{\tau}^N(H)\right|>\epsilon\right)=-\infty
\end{equation}
for any $\epsilon>0$ and $H\in C[0, 1]$, where $\mathcal{T}$ is the set of stopping times of $\{X_t^N\}_{t\geq 0}$ bounded by $T$.

We first check Equation \eqref{equ 5.5 EtightC1}. For any $H\in C[0, 1]$,
\[
\sup_{0\leq t\leq T}|\mu_t^N(H)|\leq \frac{\|H\|_{\infty}}{N}\sum_{i=1}^NX_0^N(i).
\]
Then, by Assumption (A) and Chebyshev's inequality,
\[
P\left(\sup_{0\leq t\leq T}\left|\mu_t^N(H)\right|>M\right)\leq e^{-N\left(M-\int_0^1\phi(x)dx\left(e^{\|H\|_{\infty}}-1\right)+o(1)\right)}
\]
and hence Equation \eqref{equ 5.5 EtightC1} holds.

Now we check Equation \eqref{equ 5.6 EtightC2}. For any $M\geq1$, we use $B_M^N$ to denote the event $\left\{\frac{1}{N}\sum_{N=1}^NX_0^N(i)\leq M\right\}$.
For $c>0$, by Equation \eqref{equ 5.1 Lambda},
\[
\frac{\Lambda_{t+\tau}^N(cH)}{\Lambda_{\tau}^N(cH)}=\exp\left\{cN\left(\mu^N_{t+\tau}(H)-\mu^N_{\tau}(H)-\int_{\tau}^{t+\tau}\mu_s^N(\mathcal{B}^NH)ds\right)\right\}
\]
and hence
\begin{align*}
& \left\{\sup_{0\leq t\leq \delta}\left(\mu_{t+\tau}^N(H)-\mu_{\tau}^N(H)\right)>\epsilon\right\}\bigcap B^N_M\\
&\subseteq \left\{\sup_{0\leq t\leq \delta}\frac{\Lambda_{t+\tau}^N(cH)}{\Lambda_{\tau}^N(cH)}\geq \exp\left\{cN\epsilon-cNM\delta\|\lambda\|_{\infty}(e^{2\|H\|_{\infty}}-1)\right\}\right\}.
\end{align*}
Then, since $\left\{\frac{\Lambda_{t+\tau}^N(cH)}{\Lambda_{\tau}^N(cH)}\right\}_{t\geq 0}$ is a martingale according to Lemma \ref{lemma LambdaisMartingale}, by Doob's inequality,
\[
\frac{1}{N}\log P\left(\sup_{0\leq t\leq \delta}\left(\mu_{t+\tau}^N(H)-\mu_{\tau}^N(H)\right)>\epsilon, B^N_M\right)\leq -c\epsilon+cM\delta\|\lambda\|_{\infty}(e^{2\|H\|_{\infty}}-1)
\]
and hence
\[
\limsup_{\delta\rightarrow0}\limsup_{N\rightarrow+\infty}\frac{1}{N}\log\sup_{\tau\in \mathcal{T}}P\left(\sup_{0\leq t\leq \delta}\left(\mu_{t+\tau}^N(H)-\mu_{\tau}^N(H)\right)>\epsilon, B_M^N\right)\leq -c\epsilon.
\]
Since $c$ is arbitrary,
\[
\limsup_{\delta\rightarrow0}\limsup_{N\rightarrow+\infty}\frac{1}{N}\log\sup_{\tau\in \mathcal{T}}P\left(\sup_{0\leq t\leq \delta}\left(\mu_{t+\tau}^N(H)-\mu_{\tau}^N(H)\right)>\epsilon, B_M^N\right)=-\infty
\]
and hence
\begin{align*}
&\limsup_{\delta\rightarrow0}\limsup_{N\rightarrow+\infty}\frac{1}{N}\log\sup_{\tau\in \mathcal{T}}P\left(\sup_{0\leq t\leq \delta}\left(\mu_{t+\tau}^N(H)-\mu_{\tau}^N(H)\right)>\epsilon\right)\\
&\leq \limsup_{N\rightarrow +\infty}\frac{1}{N}\log P\left(\left(B_M^N\right)^c\right).
\end{align*}
By Assumption (A) and Chebyshev's inequality,
\[
\limsup_{N\rightarrow +\infty}\frac{1}{N}\log P\left(\left(B_M^N\right)^c\right) \leq -M+(e-1)\int_0^1\phi(x)dx.
\]
Since $M$ is arbitrary, we have
\begin{equation}\label{equ 5.7}
\limsup_{\delta\rightarrow0}\limsup_{N\rightarrow+\infty}\frac{1}{N}\log\sup_{\tau\in \mathcal{T}}P\left(\sup_{0\leq t\leq \delta}\left(\mu_{t+\tau}^N(H)-\mu_{\tau}^N(H)\right)>\epsilon\right)=-\infty.
\end{equation}
For $c>0$, since $\left\{\frac{\Lambda_{t+\tau}^N(-cH)}{\Lambda_{\tau}^N(-cH)}\right\}_{t\geq 0}$ is also a martingale, an analysis similar with that leading to Equation \eqref{equ 5.7} shows that
\begin{equation}\label{equ 5.8}
\limsup_{\delta\rightarrow0}\limsup_{N\rightarrow+\infty}\frac{1}{N}\log\sup_{\tau\in \mathcal{T}}P\left(\inf_{0\leq t\leq \delta}\left(\mu_{t+\tau}^N(H)-\mu_{\tau}^N(H)\right)<-\epsilon\right)=-\infty.
\end{equation}
Equation \eqref{equ 5.6 EtightC2} follows from Equations \eqref{equ 5.7} and \eqref{equ 5.8} while the proof is complete.

\qed

\section{Proof of Equation \eqref{equ LDP lower}}\label{section six}

In this section, we give the proof of Equation \eqref{equ LDP lower}. For simplicity, in proofs we write a random variable $\varrho_N$ as $o_{exp}(\frac{1}{N})$ when $\{\varrho_N\}_{N\geq 1}$ satisfies
\[
\limsup_{N\rightarrow+\infty}\frac{1}{N}\log P\left(|\varrho_N|>\epsilon\right)=-\infty
\]
for any $\epsilon>0$. We first give a more clear expression of $I_{ini}(\mu_0)$ and $I_{dyn}(\mu)$ for $\mu\in D_0$.

\begin{lemma}\label{lemma 6.1}
Under Assumption (B), for any $\mu\in D_0$,
\begin{equation*}
I_{ini}(\mu_0)=\int_0^1\psi(0, x)\log \psi(0, x)-\psi(0, x)\log\phi(x)+\phi(x)-\psi(0, x)dx
\end{equation*}
where $\psi(0, \cdot)(x)=\frac{d\mu_0}{dx}$ and there exists $G\in C\left([0, T]\times[0, 1]\right)$ such that
\begin{equation}\label{equ 6.1 definition of G}
\partial_s\psi(s, x)=\int_0^1\psi(s, y)\lambda(y, x)e^{G_s(x)-G_s(y)}-\psi(s, x)\lambda(x, y)e^{G_s(y)-G_s(x)}dy
\end{equation}
and
\begin{align*}
I_{dyn}(\mu)=&\int_0^T\int_0^1\partial_s\psi_s(x)G_s(x)dsdx\\
&-\int_0^T\int_0^1\int_0^1\psi_s(x)\lambda(x, y)\left(e^{G_s(y)-G_s(x)}-1\right)dsdxdy,
\end{align*}
where $\psi(s, \cdot)(x)=\frac{d\mu_s}{dx}$.
\end{lemma}

\proof

For the first part of Lemma \ref{lemma 6.1},
\[
I_{ini}(\mu_0)=\int_0^1\psi(0, x)\log \psi(0, x)-\psi(0, x)\log\phi(x)+\phi(x)-\psi(0, x)dx
\]
follows from the fact that
\[
\Xi(\theta)=\psi(0, x)\theta-\phi(x)(e^\theta-1)
\]
gets maximum $\psi(0, x)\log\psi(0, x)-\psi(0, x)\log\phi(x)-\psi(0, x)+\phi(x)$ at $\theta=\log \frac{\psi(0, x)}{\phi(x)}$.

For the second part, if $\psi\equiv 0$, then it is easy to check that this part holds with $G\equiv 1$. So we only need to deal with the case where $\psi\not\equiv 0$. Since $\int_0^1\partial_s\psi(s, x) dx=0$, $\psi\not\equiv 0$ implies that $\psi(s, \cdot)\not\equiv 0$ for any $0\leq s\leq T$. Then, for any $0\leq s\leq T$, it is easy to check that there exists a unique $C_s>0$ such that
\[
2C_s=\int_0^1\sqrt{\left(\frac{\partial}{\partial s}\psi(s, y)\right)^2+4\psi(s,y)\lambda_1(y)\lambda_2(y)C_s} dy
\]
according to the fact that
\[
\varpi(x)=2x-\int_0^1\sqrt{\left(\frac{\partial}{\partial s}\psi(s, y)\right)^2+4\psi(s,y)\lambda_1(y)\lambda_2(y)x} dy
\]
is a convex function with respect to $x$ while $\varpi(0)\leq 0, \varpi(+\infty)=+\infty$. Let
\[
G_s(x)=\log\left(\frac{\partial_s\psi(s, x)+\sqrt{\left(\frac{\partial}{\partial s}\psi(s, x)\right)^2+4\psi(s,x)\lambda_1(x)\lambda_2(x)C_s}}{2C_s\lambda_2(x)}\right),
\]
then, by direct calculation, it is easy to check that $G$ satisfies Equation \eqref{equ 6.1 definition of G}.

For $\mu\in D_0$, according to the formula of integration by parts and the fact that $C^{1, 0}\left([0, T]\times [0, 1]\right)$ is dense in $C\left([0, T]\times [0, 1]\right)$,
\begin{align*}
&I_{dyn}(\mu)=\sup\Bigg\{\int_0^1\int_0^T\partial_s\psi_s(x)H_s(x)dxds\\
&-\int_0^T\int_0^1\int_0^1\psi(s, x)\lambda(x, y)\left(e^{H_s(y)-H_s(x)}-1\right)dsdxdy:~H\in C\left([0, T]\times[0, 1]\right)\Bigg\}.
\end{align*}
For simplicity, we use $\vartheta(H)$ to denote
\[
\int_0^1\int_0^T\partial_s\psi_s(x)H_s(x)dxds-\int_0^T\int_0^1\int_0^1\psi(s, x)\lambda(x, y)\left(e^{H_s(y)-H_s(x)}-1\right)dsdxdy.
\]
For any $\epsilon>0$ and any $H\in C\left([0, T]\times[0, 1]\right)$, let $k_\epsilon(G, H)=\vartheta(G+\epsilon(H-G))$, then it is easy to check that $k_\epsilon$ is a concave function with respect to $\epsilon$ and $\frac{dk_\epsilon}{d\epsilon}=0$ when $\epsilon=0$ since $G$ satisfies Equation \eqref{equ 6.1 definition of G}. Hence, $\vartheta(G)=k_0\geq k_1=\vartheta(H)$ and
\[
I_{dyn}(\mu)=\sup_{H\in C\left([0, T]\times[0, 1]\right)}\vartheta(H)=\vartheta(G).
\]

\qed

As a second step, we need to give the hydrodynamic limit of $\{X_t^N\}_{t\geq 0}$ under a transformed measure with the exponential martingale introduced in Section \ref{section five} as R-N derivative with respect to the original probability measure $P$.  To give the precise statement of our result, we introduce some notations. For any positive $\gamma\in C[0, 1]$, we use $P_{\gamma}$ to denote the probability measure of our process $\{X_t^N\}_{t\geq 0}$ with initial condition where $\{X_0^N(i)\}_{1\leq i\leq N}$ are independent and $X_0^N(i)$ follows Poisson distribution with mean $\gamma(\frac{i}{N})$ for each $1\leq i\leq N$. For any $G\in C^{1,0}\left([0, T]\times [0, 1]\right)$, we use $\widehat{P}_\gamma^G$ to denote the probability measure such that
\[
\frac{d\widehat{P}_\gamma^G}{dP_\gamma}=\Lambda_T^N(G).
\]
For positive $\gamma\in C[0, 1]$ and $G\in C^{1,0}\left([0, T]\times [0, 1]\right)$, according to the theory of ordinary differential equations on Banach spaces, it is easy to check that there exists an unique $\gamma^G\in C^{1,0}\left([0, T]\times [0, 1]\right)$ such that
\begin{equation}\label{equ 6.2 definition of gammaG}
\begin{cases}
&\partial_s\gamma^G(s, x)=\int_0^1\gamma^G(s, y)\lambda(y, x)e^{G_s(x)-G_s(y)}-\gamma^G(s, x)\lambda(x, y)e^{G_s(y)-G_s(x)}dy \\
&\text{~\quad\quad\quad}\text{~for all~}0\leq x\leq 1 \text{~and~} 0\leq s\leq T,\\
&\gamma^G(0, x)=\gamma(x)  \text{~for all~}0\leq x\leq 1.
\end{cases}
\end{equation}
Then, we have the following lemma.

\begin{lemma}\label{lemma 6.2 hydrodynamic}
For positive $\gamma\in C[0, 1]$, $\{\mu^N\}_{N\geq 1}$ converges in $\widehat{P}_{\gamma}^G$-probability to $\mu^{\gamma, G}$, where
\begin{equation}\label{equ 6.9 definition of muGammaG}
\mu^{\gamma, G}_t(dx)=\gamma^G(t, x)dx
\end{equation}
for $0\leq t\leq T$.
\end{lemma}
To prove Lemma \ref{lemma 6.2 hydrodynamic}, we need the following lemma.
\begin{lemma}\label{lemma 6.3 error}
For given $G\in C^{1,0}\left([0, T]\times[0, 1]\right)$ and positive $\gamma\in C[0, 1]$, under $\widehat{P}_\gamma^G$,
\begin{equation}\label{equ 6.3 error exponentially small}
\sup_{H\in C[0, 1], \|H\|_{\infty}=1}\sum_{0\leq t\leq T}\left(\mu^N_t(H)-\mu_{t-}^N(H)\right)^2=o_{exp}(\frac{1}{N}).
\end{equation}
\end{lemma}

\proof[Proof of Lemma \ref{lemma 6.3 error}]

For any $H\in C[0, 1]$ with $\|H\|_{\infty}=1$, on each jumping moment $s$ of $\{X_t^N\}_{t\geq 0}$,
\[
|\mu_s^N(H)-\mu_{s-}^N(H)|^2\leq \left(\frac{2\|H\|_{\infty}}{N}\right)^2=\frac{4}{N^2}.
\]
For simplicity, we use $\varsigma^N$ to denote $\sup_{H\in C[0, 1], \|H\|_{\infty}=1}\sum_{0\leq t\leq T}\left(\mu^N_t(H)-\mu_{t-}^N(H)\right)^2$. For $M\geq 1$, since each ball jumps at rate at most $\|\lambda\|_{\infty}$ under $P_\gamma$, condition on
\[
\frac{1}{N}\sum_{i=1}^NX_0^N(i)\leq M,
\]
$\varsigma^N$ is stochastically dominated from above by $\frac{4}{N^2}Y(NM\|\lambda\|_\infty T)$ under $P_\gamma$, where $\{Y(t)\}_{t\geq 0}$ is the Poisson process with rate $1$ defined as in Section \ref{section three}. Therefore, by Chebyshev's inequality,
\[
P_\gamma\left(\varsigma^N\geq \epsilon, \frac{1}{N}\sum_{i=1}^NX_0^N(i)\leq M\right)\leq e^{-\frac{\epsilon N^2}{4}}e^{(e-1)\|\lambda\|_{\infty}MNT}
\]
and
\[
\limsup_{N\rightarrow+\infty}\frac{1}{N}\log P_\gamma\left(\varsigma^N\geq \epsilon, \frac{1}{N}\sum_{i=1}^NX_0^N(i)\leq M\right)=-\infty.
\]
Therefore,
\[
\limsup_{N\rightarrow+\infty}\frac{1}{N}\log P_\gamma\left(\varsigma^N\geq \epsilon\right)
\leq \limsup_{N\rightarrow+\infty}\frac{1}{N}\log P_\gamma\left(\frac{1}{N}\sum_{i=1}^NX_0^N(i)\geq M\right).
\]
Under $P_\gamma$,
\[
P_\gamma\left(\frac{1}{N}\sum_{i=1}^NX_0^N(i)\geq M\right)\leq e^{-MN}e^{(e-1)\sum_{i=1}^N\gamma(\frac{i}{N})}
\]
and hence
\[
\limsup_{N\rightarrow+\infty}\frac{1}{N}\log P_\gamma\left(\varsigma^N\geq \epsilon\right)
\leq -M+(e-1)\int_0^1\gamma(x)dx.
\]
Since $M$ is arbitrary, Equation \eqref{equ 6.3 error exponentially small} holds under $P_\gamma$. By Equation \eqref{equ 5.2 two}, it is easy to check that
\[
{\rm E}_{\gamma}\left(\left(\Lambda_T^N(G)\right)^2\right)=e^{NO(1)}.
\]
Then, Equation \eqref{equ 6.3 error exponentially small} holds under $\widehat{P}_\gamma^G$ according to H\"{o}lder's inequality.

\qed

Now we give the proof of Lemma \ref{lemma 6.2 hydrodynamic}. The proof follows a strategy similar with those introduced in References \cite{Xue2019+, XueZhao2020+b}, where a generalized version of Girsanov's theorem introduced in \cite{Schuppen1974} will be utilized.

\proof[Proof of Lemma \ref{lemma 6.2 hydrodynamic}]

According to the definition of $\mu^{\gamma, G}$ and the formula of integration by parts, for any $0\leq t\leq T$ and $H\in C[0, 1]$ with $\|H\|_\infty=1$,
\begin{equation}\label{equ 6.4 integration by parts}
\mu^{\gamma, G}_t(H)=\int_0^1\gamma(x)H(x)dx+\int_0^t\mu_s^{\gamma, G}\left(\upsilon^G_sH\right)ds,
\end{equation}
where
\[
(\upsilon_s^GH)(x)=\int_0^1\lambda(x, y)e^{G_s(y)-G_s(x)}\left(H(y)-H(x)\right)dy
\]
for any $x\in C[0, 1]$.

For any $H\in C[0, 1]$ with $\|H\|_{\infty}=1$, we define $\Theta_H^N\in C\left(X\right)$ as
\[
\Theta_H^N(x)=\frac{1}{N}\sum_{i=1}^Nx(i)H(\frac{i}{N}),
\]
then $\Theta_H^N(X_t^N)=\mu_t^N(H)$ and, as we have defined in Section \ref{section five},
\[
\mathcal{U}_t^N(\Theta_H^N)=\mu_t^N(H)-\mu_0^N(H)-\int_0^t\mathcal{L}_N\mu_s^N(H)ds
\]
is a martingale under $P_\gamma$. Let $f_G^N$ be defined as in Section \ref{section five}, then we define
\[
\widetilde{\mathcal{U}}_t^N(f_G^N)=\int_0^t\frac{1}{f_G^N(s-, X_{s-}^N)}d\mathcal{U}_s^N(f_G^N)
\]
for $0\leq t\leq T$. Then, by Equation \eqref{equ 5.2 dLambda},
\begin{equation}\label{equ 6.5 dLambda}
d\Lambda_t^N(G)=\Lambda_{t-}^N(G)d\widetilde{\mathcal{U}}_t^N(f_G^N).
\end{equation}
By Equation \eqref{equ 6.5 dLambda}, Lemma \ref{lemma LambdaisMartingale} and Theorem 3.2 of \cite{Schuppen1974},
\[
\widehat{\mathcal{U}}_t^N(\Theta_H^N):=\mathcal{U}_t^N(\Theta_H^N)-\langle\mathcal{U}^N(\Theta_H^N), \widetilde{\mathcal{U}}^N(f_G^N)\rangle_t
\]
is a martingale under $\widehat{P}_\gamma^G$ and
\[
\left[\widehat{\mathcal{U}}^N(\Theta_H^N), \widehat{\mathcal{U}}^N(\Theta_H^N)\right]_t=[\mathcal{U}^N(\Theta_H^N), \mathcal{U}^N(\Theta_H^N)]_t
=\sum_{0\leq s\leq t}\left(\mu_s^N(H)-\mu^N_{s-}(H)\right)^2
\]
under both $P_\gamma$ and $\widehat{P}_\gamma^G$. By Lemma \ref{lemma 6.3 error} and Doob's inequality,
\begin{equation}\label{equ 6.6 martinagle is small}
\sup_{0\leq t\leq T, \atop H\in C[0, 1], \|H\|_\infty=1}\left|\widehat{\mathcal{U}}_t^N(\Theta_H^N)\right|=o_p(1)
\end{equation}
under $\widehat{P}^G_\gamma$. According to the definition of $\widetilde{\mathcal{U}}_t^N(f_G^N)$,
\[
d\langle\mathcal{U}^N(\Theta_H^N), \widetilde{\mathcal{U}}^N(f_G^N)\rangle_t=\frac{1}{f_G^N(t-, X_{t-}^N)}d\langle\mathcal{U}^N(\Theta_H^N), \mathcal{U}^N(f_G^N)\rangle_t.
\]
By Dynkin's martingale formula,
\[
d\langle\mathcal{U}^N(\Theta_H^N), \mathcal{U}^N(f_G^N)\rangle_t=\left(\mathcal{L}_N(f_G^N\Theta_H^N)-\Theta_H^N\mathcal{L}_Nf_G^N-f_G^N\mathcal{L}_N\Theta_H^N\right)(t, X_t)dt.
\]
By the definition of $\mathcal{L}_N$ and direct calculations, it is not difficult to check that
\begin{align*}
&\left(\mathcal{L}_N(f_G^N\Theta_H^N)-\Theta_H^N\mathcal{L}_Nf_G^N-f_G^N\mathcal{L}_N\Theta_H^N\right)(t, X_t^N)\\
&=f_G^N(t, X_t^N)\left(\mu_t^N(\upsilon_t^{G, N}H)-\mathcal{L}_N\mu_t^N(H)\right),
\end{align*}
where
\[
(\upsilon_t^{G, N}H)(x)=\sum_{j=1}^N\lambda\left(x, \frac{j}{N}\right)e^{G_t(\frac{j}{N})-G_t(x)}\left(H(\frac{j}{N})-H(x)\right)
\]
for any $x\in [0, 1]$. As a result, under $\widehat{P}_\gamma^G$,
\begin{align*}
\mu_t^N(H)&=\mu_0^N(H)+\int_0^t\mu^N_s(\upsilon_s^{G, N}H)ds+o_p(1)\\
&=\int_0^1\gamma(x)H(x)dx+\int_0^t\mu^N_s(\upsilon_s^{G, N}H)ds+o_p(1).
\end{align*}
Note that in the above equation we utilize the fact that $\mu_0^N(H)=\int_0^1\gamma(x)H(x)dx+o_p(1)$ under $\widehat{P}_\gamma^G$ and this $o_p(1)$ can be chosen uniformly for all $H\in C[0, 1]$ with $\|H\|_\infty=1$.

According to an analysis similar with that given in the proof of Lemma \ref{lemma 5.1 replacement}, it is not difficult to check that
\[
\sup_{0\leq t\leq T,\atop H\in C[0, 1], \|H\|_\infty=1}\left|\mu_t^N(\upsilon_t^{G, N}H)-\mu_t^N(\upsilon_t^GH)\right|=o_{exp}(\frac{1}{N})
\]
under $P_\gamma$. Then, according to H\"{o}lder's inequality,
\[
\sup_{0\leq t\leq T,\atop H\in C[0, 1], \|H\|_\infty=1}\left|\mu_t^N(\upsilon_t^{G, N}H)-\mu_t^N(\upsilon_t^GH)\right|=o_{exp}(\frac{1}{N})=o_p(1)
\]
under $\widehat{P}_\gamma^G$. As a result, for $0\leq t\leq T$ and $H\in C[0, 1]$ with $\|H\|_\infty=1$,
\begin{equation}\label{equ 6.7}
\mu_t^N(H)=\int_0^1\gamma(x)H(x)dx+\int_0^t\mu_s^N(\upsilon_s^GH)ds+o_p(1),
\end{equation}
where $o_p(1)$ can be chosen uniformly for all $H\in C[0, 1]$ with $\|H\|_\infty$. According to Equations \eqref{equ 6.4 integration by parts}, \eqref{equ 6.7} and the fact that $\|\upsilon_s^GH\|_\infty\leq 2\|H\|_\infty\|\lambda\|_{\infty}e^{2\|G\|_\infty}$,
\begin{align*}
&\sup_{H\in C[0, 1], \|H\|=1}\left|\mu_t^{\gamma, G}(H)-\mu_t^N(H)\right|\\
&\leq o_p(1)+2\|\lambda\|_{\infty}e^{2\|G\|_\infty}
\int_0^t \sup_{H\in C[0, 1], \|H\|=1}\left|\mu_s^{\gamma, G}(H)-\mu_s^N(H)\right|ds
\end{align*}
for any $0\leq t\leq T$ under $\widehat{P}_\gamma^G$. Then, by Grown-wall's inequality,
\begin{equation}\label{equ 6.8}
\sup_{H\in C[0, 1], \|H\|=1}\left|\mu_t^{\gamma, G}(H)-\mu_t^N(H)\right|\leq o_p(1)\exp\left\{2t\|\lambda\|_{\infty}e^{2\|G\|_\infty}\right\}
\end{equation}
under $\widehat{P}_\gamma^G$. Lemma \ref{lemma 6.2 hydrodynamic} follows from Equation \eqref{equ 6.8} directly.

\qed

At last we give the proof of Equation \eqref{equ LDP lower}.

\proof[Proof of Equation \eqref{equ LDP lower}]

For any $\epsilon>0$, there exists $\mu^\epsilon\in D_0\bigcap O$ such that
\[
I_{ini}(\mu^\epsilon_0)+I_{dyn}(\mu^\epsilon)\leq \inf_{\mu\in D_0\bigcap O}\left(I_{ini}(\mu_0)+I_{dyn}(\mu)\right)+\epsilon.
\]
By Lemma \ref{lemma 6.1}, there exists $G^\epsilon\in C\left([0, T]\times [0, 1]\right)$ such that
\[
\partial_s\psi^\epsilon(s, x)=\int_0^1\psi^\epsilon(s, y)\lambda(y, x)e^{G^\epsilon_s(x)-G^\epsilon_s(y)}-\psi^\epsilon(s, x)\lambda(x, y)e^{G^\epsilon_s(y)-G^\epsilon_s(x)}dy
\]
and
\begin{align*}
I_{dyn}(\mu^\epsilon)=&\int_0^T\int_0^1\partial_s\psi^\epsilon_s(x)G^\epsilon_s(x)dsdx\\
&-\int_0^T\int_0^1\int_0^1\psi^\epsilon_s(x)\lambda(x, y)\left(e^{G^\epsilon_s(y)-G^\epsilon_s(x)}-1\right)dsdxdy\\
=&\mu^\epsilon_T(G^\epsilon_T)-\mu^\epsilon_0(G^\epsilon_0)-\int_0^T\mu^\epsilon_s\left((\partial_s+\mathcal{B})G^\epsilon_s\right)ds
\end{align*}
while
\[
I_{ini}(\mu^\epsilon_0)=\int_0^1\psi^\epsilon(0, x)\log \psi^\epsilon(0, x)-\psi^\epsilon(0, x)\log\phi(x)+\phi(x)-\psi^\epsilon(0, x)dx,
\]
where $\psi^\epsilon_s=\frac{d\mu^\epsilon_s}{dx}$.

Since $C^{1,0}\left([0, T]\times [0, 1]\right)$ is dense in $C\left([0, T]\times [0, 1]\right)$ under the $l_\infty$ norm, there exists $C^{1,0}\left([0, T]\times [0, 1]\right)$-valued sequence $\{G^n\}_{n\geq 1}$ such that $\|G^n-G^\epsilon\|_\infty\rightarrow 0$ as $n\rightarrow+\infty$. For each $n\geq 1$ and $x\in [0, 1]$, we define $\gamma^n(x)=\psi^\epsilon(0, x)+\frac{1}{n}$ and then $\gamma^n$ is strictly positive. We use $\zeta^n$ to denote $\mu^{\gamma^n, G^n}$ defined as in Equation \eqref{equ 6.9 definition of muGammaG}. Then, by Lemma \ref{lemma 6.1},
\begin{align}\label{equ 6.11 IdynMun}
I_{dyn}(\zeta^n)=&\int_0^T\int_0^1\partial_s\psi^n_s(x)G^n_s(x)dsdx \notag\\
&-\int_0^T\int_0^1\int_0^1\psi^n_s(x)\lambda(x, y)\left(e^{G^n_s(y)-G^n_s(x)}-1\right)dsdxdy\notag\\
=&\zeta^n_T(G^n_T)-\zeta^n_0(G^n_0)-\int_0^T\zeta^n_s\left((\partial_s+\mathcal{B})G^n_s\right)ds
\end{align}
and
\[
I_{ini}(\zeta^n_0)=\int_0^1\gamma^n(x)\log \gamma^n(x)-\gamma^n(x)\log\phi(x)+\phi(x)-\gamma^n(0, x)dx,
\]
where $\psi_s^n=\frac{d\zeta^n_s}{dx}, \psi_0^n=\frac{d\zeta^n_0}{dx}=\gamma^n$ while
\[
\partial_s\psi^n(s, x)=\int_0^1\psi^n(s, y)\lambda(y, x)e^{G^n_s(x)-G^n_s(y)}-\psi^n(s, x)\lambda(x, y)e^{G^n_s(y)-G^n_s(x)}dy.
\]
Furthermore, according to Grownwall's inequality, it is not difficult to check that $\zeta^n\rightarrow \mu^\epsilon$ in $D\big([0, T], \mathcal{S}\big)$ and $I_{dyn}(\zeta^n)\rightarrow I_{dyn}(\mu^\epsilon)$, $I_{ini}(\zeta^n_0)\rightarrow I_{ini}(\mu^\epsilon_0)$. Since $O$ is open while $\mu^\epsilon\in O$, there exists integer $m=m(\epsilon)\geq 1$ such that
\[
I_{ini}(\zeta^m_0)+I_{dyn}(\zeta^m)\leq I_{ini}(\mu^\epsilon_0)+I_{dyn}(\mu^\epsilon)+\epsilon \text{~and~} \mu^m\in O.
\]
For $N\geq 1$, we use $W_N$ to denote the subset
\[
\left\{\mu\in D\big([0, T], \mathcal{S}\big):~\sup_{0\leq s\leq T}\left|\mu_s(\mathcal{B}G_s^m)-\mu_s(\mathcal{B}^NG_s^m)\right|<\epsilon\right\}
\]
of $D\big([0, T], \mathcal{S}\big)$ while use $Q$ to denote the subset
\begin{align*}
\Bigg\{\mu\in D\big([0, T], \mathcal{S}\big): \text{~}&\Bigg|(\zeta^m-\mu)_T(G^m_T)\\
&-(\zeta^m-\mu)_0(G^m_0)-\int_0^T(\zeta^m-\mu)_s\left((\partial_s+\mathcal{B})G^m_s\right)ds\Bigg|<\epsilon\Bigg\}
\end{align*}
of $D\big([0, T], \mathcal{S}\big)$. According to Lemma \ref{lemma 5.1 replacement} and H\"{o}lder's inequality, it is easy to check that
\begin{equation}\label{equ 6.10 WN high probability}
\lim_{N\rightarrow+\infty}\widehat{P}_{\gamma^m}^{G^m}\left(\mu^N\in W_N\right)=\lim_{N\rightarrow+\infty}P_{\gamma^m}\left(\mu^N\in W_N\right)=1.
\end{equation}
By Lemma \ref{lemma 6.2 hydrodynamic}, $\mu^N\rightarrow \mu^{\gamma^m, G^m}=\zeta^m\in O$ in $\widehat{P}_{\gamma^m}^{G^m}$-probability and hence
\[
\lim_{N\rightarrow+\infty}\widehat{P}_{\gamma^m}^{G^m}\left(\mu^N\in Q\bigcap O\right)=1.
\]
Then, by Equation \eqref{equ 6.10 WN high probability},
\begin{equation}\label{equ 6.11 WNandQ high probability}
\lim_{N\rightarrow+\infty}\widehat{P}_{\gamma^m}^{G^m}\left(\mu^N\in W_N\bigcap Q\bigcap O\right)=1.
\end{equation}
On the event $\{\mu^N\in W_N\bigcap Q\}$, by Equations \eqref{equ 5.1 Lambda} and \eqref{equ 6.11 IdynMun},
\[
\Lambda_T^N(G^m)\leq \exp\left\{N\left(I_{dyn}(\zeta^m)+(T+1)\epsilon\right)\right\}
\]
and hence
\[
\frac{dP_{\gamma^m}}{d\widehat{P}^{G_m}_{\gamma^m}}\geq \exp\left\{-N\left(I_{dyn}(\zeta^m)+(T+1)\epsilon\right)\right\}.
\]
According to the definition of $P_{\gamma^m}$ and Assumption (A),
\begin{align*}
\frac{dP}{dP_{\gamma^m}}&=\prod_{i=1}^n\left[e^{-\phi(\frac{i}{n})+\gamma^m(\frac{i}{n})}\left(\frac{\phi(\frac{i}{N})}{\gamma^m(\frac{i}{N})}\right)^{X_0^N(i)}\right]\\
&=\exp\left\{-N\left[\mu_0^N(\log\gamma^m-\log\phi)+\int_0^1\phi(x)dx-\int_0^1\gamma^m(x)dx+o(1)\right]\right\}.
\end{align*}
Under $\widehat{P}^{G^m}_{\gamma^m}$, since $\zeta_0^m(dx)=\gamma^m(x)dx$, it is easy to check that $\mu_0^N\rightarrow \zeta_0^m$. As a result, let $Q_1$ be the subset
\[
\left\{\mu\in D\big([0, T], \mathcal{S}\big):~\left|\mu_0(\log\phi-\log\gamma^m)-\zeta_0^m(\log\phi-\log\gamma^m)\right|<\epsilon\right\}
\]
of $D\big([0, T], \mathcal{S}\big)$, then
\[
\lim_{N\rightarrow+\infty}\widehat{P}_{\gamma^m}^{G^m}\left(\mu^N\in Q_1\right)=1
\]
and hence
\begin{equation}\label{equ 6.13}
\lim_{N\rightarrow+\infty}\widehat{P}_{\gamma^m}^{G^m}\left(\mu^N\in O\bigcap W_N\bigcap Q\bigcap Q_1\right)=1
\end{equation}
by Equation \eqref{equ 6.11 WNandQ high probability}. On $Q_1$, by Lemma \ref{lemma 6.1},
\[
\frac{dP}{dP_{\gamma^m}}\geq \exp\left\{-N\left[I_{ini}(\zeta_0^m)+\epsilon+o(1)\right]\right\}.
\]
Hence, on $W_N\bigcap Q\bigcap Q_1$,
\[
\frac{dP}{d\widehat{P}_{\gamma^m}^{G^m}}=\frac{dP}{dP_{\gamma^m}}\frac{dP_{\gamma^m}}{d\widehat{P}^{G_m}_{\gamma^m}}
\geq \exp\left\{-N\left[I_{dyn}(\zeta^m)+I_{ini}(\zeta_0^m)+(T+2)\epsilon+o(1)\right]\right\}.
\]
As a result, by Equation \eqref{equ 6.13},
\begin{align*}
P\left(\mu^N\in O\right)&\geq P\left(\mu^N\in O\bigcap W_N\bigcap Q\bigcap Q_1\right) \\
&={\rm E}_{\gamma^m}^{G^m}\left(\frac{dP}{d\widehat{P}_{\gamma^m}^{G^m}}1_{_{O\bigcap W_N\bigcap Q\bigcap Q_1}}\right)\\
&\geq \exp\left\{-N\left[I_{dyn}(\zeta^m)+I_{ini}(\zeta_0^m)+(T+2)\epsilon+o(1)\right]\right\}\left(1+o(1)\right).
\end{align*}
Therefore,
\begin{align*}
\liminf_{N\rightarrow+\infty}\frac{1}{N}\log P\left(\mu^N\in O\right)&\geq -\left(I_{dyn}(\zeta^m)+I_{ini}(\zeta_0^m)\right)-(T+2)\epsilon\\
&\geq -\left(I_{dyn}(\mu^\epsilon)+I_{ini}(\mu_0^\epsilon)\right)-(T+3)\epsilon\\
&\geq -\inf_{\mu\in D_0\bigcap O}\left(I_{ini}(\mu_0)+I_{dyn}(\mu)\right)-(T+4)\epsilon.
\end{align*}
Since $\epsilon$ is arbitrary, Equation \eqref{equ LDP lower} holds.

\qed

\appendix{}
\section{Appendix}
\subsection{Proof of Lemma \ref{Existen and uniqueness}}

\proof[Proof of Lemma \ref{Existen and uniqueness}]

We first construct a weak solution $\mu$ of Equation \eqref{equ 2.1 ODE on banach} to show the existence. We use $\|f\|_\infty$ to denote the $l_\infty$-norm of $f$ for any $f\in C[0, 1]$ and $\|\lambda\|_{\infty}$ to denote $\sup_{0\leq x, y\leq 1}\lambda(x, y)$. For each $f\in C[0, 1]$, we define
\[
\left(P_3f\right)(x)=\int_0^1\lambda(y, x)f(y)dy
\]
for any $x\in [0, 1]$. Since $\|P_2f\|_\infty, \|P_3f\|_\infty\leq \|\lambda\|_{\infty}\|f\|_\infty$ for any $f\in C[0, 1]$, the linear $C[0,1]$-valued ODE
\[
\begin{cases}
\frac{d}{dt}\rho_t=(P_3-P_2)\rho_t, &\text{~}0\leq t\leq T,\\
\rho_0(x)=\phi(x), &\text{~}0\leq x\leq 1,
\end{cases}
\]
satisfies Lipschitz' condition and hence has the unique solution
\[
\rho_t=e^{t(P_3-P_2)}\phi=\sum_{k=0}^{+\infty}\frac{t^k\left(P_3-P_2\right)^k}{k!}\phi.
\]
Let $\mu_t(dx)=\rho_t(x)dx$ for any $0\leq t\leq T$, then $\mu=\{\mu_t\}_{0\leq t\leq T}$ is a weak solution of Equation \eqref{equ 2.1 ODE on banach} since
\[
\int_0^1f(x)\left(P_3\rho_t\right)(x)dx=\int_0^1 \left(P_1f\right)(x)\rho_t(x)dx.
\]

Now we only need to show the uniqueness. Let $\mu^1, \mu^2$ be two weak solutions of Equation \eqref{equ 2.1 ODE on banach} and $\nu=\mu^1-\mu^2$, then
\[
\nu_t(f)=\int_0^t \nu_s\left((P_1-P_2)f\right)ds
\]
for any $0\leq t\leq T$. For any $\mathscr{A}\in \mathcal{S}$, we define $\|\mathscr{A}\|_{\infty}=\sup\left\{\mathscr{A}(f):~\|f\|_{\infty}=1\right\}$.
Then, for any $f$ satisfying $\|f\|_{\infty}=1$,
\[
|\nu_s\left((P_1-P_2)f\right)|\leq \|\nu_s\|_{\infty}\|(P_1-P_2)f\|_1\leq 2\|\nu_s\|_\infty\|\lambda\|_{\infty}.
\]
Therefore, for any $0\leq t\leq T$,
\[
\|\nu_t\|_{\infty}\leq 2\|\lambda\|_{\infty}\int_0^t \|\nu_s\|_{\infty}ds.
\]
Then, by Grown-Wall's inequality, $\|\nu_t\|_{\infty}\leq 0e^{2t\|\lambda\|_{\infty}}=0$ for any $0\leq t\leq T$ and hence $\mu^1=\mu^2$.

\qed

\quad

\textbf{Acknowledgments.}The author is grateful to Dr. Zhengyao Sun, Dr. Linjie Zhao, Prof. Minzhi Zhao and Prof. Qiang Yao for useful suggestions and comments. The author is grateful to the financial support from the National Natural Science Foundation of China with grant number 11501542.

{}
\end{document}